\documentclass[review, 3p, compress]{elsarticle}
\usepackage{hyperref}
\usepackage{amsfonts}
\usepackage{amssymb}
\usepackage{amsmath}
\usepackage{cases}
\usepackage{graphicx}
\usepackage{rotating}
\usepackage{float}
\usepackage{booktabs}
\usepackage{multirow}
\usepackage{multicol}
\usepackage{epstopdf}
\usepackage{subfigure}
\usepackage{color}
\usepackage{setspace}
\usepackage{bm}
\usepackage{doi}
\usepackage{algorithm}
\usepackage{algorithmic}
\usepackage{natbib}

\makeatletter
\@addtoreset{equation}{section}
\makeatother

\newtheorem{theorem}{Theorem}[section]
\newtheorem{lemma}{Lemma}[section]
\newtheorem{remark}{Remark}

\journal{Journal of \LaTeX\ Templates Templates}

\begin{document}

\begin{frontmatter}

\title{A limited-memory block bi-diagonal Toeplitz preconditioner
for block lower triangular Toeplitz system from time-space fractional diffusion equation}

\author[address1]{Yong-Liang Zhao\corref{correspondingauthor}}
\ead{uestc\_ylzhao@sina.com}

\author[address1]{Pei-Yong Zhu}
\ead{zpy6940@uestc.edu.cn}

\author[address2]{Xian-Ming Gu\corref{correspondingauthor}}
\cortext[correspondingauthor]{Corresponding author}
\ead{guxianming@live.cn}

\author[address1]{Xi-Le Zhao}
\ead{xlzhao122003@163.com}

\author[address3]{Jianxiong Cao}
\ead{caojianxiong2007@126.com}

\address[address1]{School of Mathematical Sciences, \\
University of Electronic Science and Technology of China, \\
Chengdu, Sichuan 611731, P.R. China}
\address[address2]{School of Economic Mathematics/Institute of Mathematics,\\
Southwestern University of Finance and Economics,\\
Chengdu, Sichuan 611130, P.R. China}
\address[address3]{School of Science, Lanzhou University of Technology, \\
Lanzhou, Gansu 730050, P.R. China}

\begin{abstract}
A block lower triangular Toeplitz system arising from
the time-space fractional diffusion equation is discussed.
For efficient solutions of such the linear system, the
preconditioned biconjugate gradient stabilized method
and the flexible general minimal residual method are exploited.
The main contribution of this paper has two aspects:
(i) A block bi-diagonal Toeplitz preconditioner is developed for the block lower triangular Toeplitz system,
whose storage is of $\mathcal{O}(N)$ with $N$ being the spatial grid number;
(ii) A new skew-circulant preconditioner is designed to accelerate
the inverse of the block bi-diagonal Toeplitz preconditioner multiplying a vector.
Numerical experiments are given to demonstrate the effectiveness of our two proposed preconditioners.
\end{abstract}


\end{frontmatter}


\section{Introduction}
\label{sec1}

In recent decades, the applications of fractional partial differential equations (FPDEs) have been interested
and recognized in numerous fields such as control systems \cite{machado2001discrete},
quantum mechanics \cite{laskin2000fractional}, stochastic dynamics \cite{gao2014mean}
and image processing \cite{bai2007fractional}.
Actually, the closed-form analytical solutions of FPDEs can be obtained in a few special cases \cite{podlubny1998fractional},
but such solutions are usually impractical.
It thus becomes imperative to study the numerical solutions of FPDEs,
and numerous reliable numerical methods have been developed
\cite{Gu2017fast, Gao2016Two, li2017galerkin, luo2016quadratic, Mao2016Efficient, hao2016finite, Liu2014A,
	cui2015compact, ccelik2012crank, lei2013circulant, zhao2014fourth,  pang2012multigrid}.
Due to the nonlocality of the fractional operators, using
the finite difference method to solve space/time-space fractional differential equations
leads to a time-stepping scheme with dense coefficient matrices.
The conventional time-stepping schemes utilizing the Gaussian elimination require the computational cost of $\mathcal{O}(N^3)$
and storage of $\mathcal{O}(N^2)$ at each time step, where $N$ is the spatial grid number.
For the purpose of optimizing the computational complexity, numerous fast algorithms
\cite{Gu2017fast, li2017galerkin, lei2013circulant, pang2012multigrid,
	Gu2014On, gu2015strang, wang2010direct, zhao2018delay} are designed.

From another point of view, if all time steps are stacked in a vector,
we will obtain an all-at-once system
or a block lower triangular system.
Ke et al. \cite{ke2015fast} combined the block forward substitution (BFS) method
with the divide-and-conquer strategy to solve the block lower triangular Toeplitz-like
with tri-diagonal blocks (BL3TB-like) system.
The complexity and storage requirement of their method are respectively
$\mathcal{O}(M N \log^2 M)$ and $\mathcal{O}(MN)$, where $M$ is the number of time steps.
Lu et al. \cite{lu2015approximate} proposed a fast approximate inversion method,
whose computational cost is of $\mathcal{O}(M N \log M)$ and storage requirement is of $\mathcal{O}(MN)$,
for the block lower triangular Toeplitz with tri-diagonal blocks (BL3TB) matrix.
The idea of this method is to approximate the coefficient matrix by the block $\epsilon$-circulant matrix,
which can be block-diagonalized by the fast Fourier transform (FFT).
Additionally, the error estimation given in \cite{lu2015approximate} shows that their method has high accuracy.
Since the sufficient condition provided in \cite{lu2015approximate} is difficult to verify in practice,
Lu et al. \cite{lu2018approximate} proposed a new sufficient condition, which is easier to check
and can be applied to several existing numerical schemes.
Huang et al. \cite{huang2017fast} combined the divide-and-conquer technique
with the circulant-and-skew-circulant representation of Toeplitz matrix inversion
for solving the nonsingular block lower triangular Toeplitz
with dense Toeplitz blocks (BLDTB) system.
Their proposed method requires a complexity within
	$\mathcal{O}\left( MN \log M \left( \log M + \log N \right) \right)$.

In this work, we mainly concentrate on fast solving the block lower triangular Toeplitz (BLTT) system
arising from time-space fractional diffusion equation (TSFDE):
\begin{equation}
	\begin{cases}
		\sideset{_0^C}{^\alpha_t}{\mathop{\mathcal{D}}} u(x,t)
		= e_{1} \sideset{_0}{^\beta_x}{\mathop{\mathcal{D}}} u(x,t)
		+ e_{2} \sideset{_x}{^\beta_L}{\mathop{\mathcal{D}}} u(x,t) + f(x,t), & 0 < t \leq T, ~0 \leq x \leq L, \\
		u(x,0) = u_{0}(x), & 0 \leq x \leq L,\\
		u(0,t) = u(L,t) = 0, & 0 \leq t \leq T,
	\end{cases}
	\label{eq1.1}
\end{equation}
where $e_1, e_2 > 0$.
The time and space fractional derivatives are introduced in Caputo and Riemann-Liouville sense \cite{podlubny1998fractional},
respectively, i.e.,
\begin{equation*}
	\sideset{_0^C}{^\alpha_t}{\mathop{\mathcal{D}}} u(x,t) =
	\frac{1}{\Gamma(1 - \alpha)} \int_{0}^{t} (t - \eta)^{-\alpha}
	\frac{\partial u(x,\eta)}{\partial \eta} d\eta, ~0 < \alpha < 1,
\end{equation*}

\begin{equation*}
	\sideset{_0}{^\beta_x}{\mathop{\mathcal{D}}} u(x,t) = \frac{1}{\Gamma(2 - \beta)} \frac{d^{2}}{dx^{2}}
	\int_{0}^{x} \frac{u(\eta,t)}{(x - \eta)^{\beta - 1}} d \eta, ~1 < \beta < 2,
\end{equation*}

\begin{equation*}
	\sideset{_x}{^\beta_L}{\mathop{\mathcal{D}}} u(x,t) =  \frac{1}{\Gamma(2 - \beta)} \frac{d^{2}}{dx^{2}}
	\int_{x}^{L} \frac{u(\eta,t)}{(\eta - x)^{\beta - 1}} d \eta, ~1 < \beta < 2,
\end{equation*}
where $\Gamma(\cdot)$ denotes the Gamma function.

In this study, we adopt the preconditioned biconjugate gradient stabilized (PBiCGSTAB) method \cite{van1992bicgstab}
and flexible generalized minimal residual (FGMRES) method
\footnote{The preconditioned sub-system is solved inexactly in each preconditioned iteration step,
	and this information just matches the characteristic of FGMRES method.
	Thus the FGMRES method is chosen in this study.} \cite{saad2003iter} to solve the BLTT system efficiently.
Therefore, the main contribution of this work can be concluded as:

(i) A block bi-diagonal Toeplitz (B2T) preconditioner, whose storage is of $\mathcal{O}(N)$,
is developed to solve the BLTT system;

(ii) A new skew-circulant preconditioner is designed to efficiently compute the inverse of
the B2T preconditioner multiplying a vector.
Furthermore, numerical experiments indicate that our skew-circulant preconditioner
is slightly better than the Strang's circulant preconditioner \cite{ng2004iterative,chan2007introduction}.

The rest of this paper is organized as follows.
In Section 2, the BLTT system is established
through the $L2$-$1_\sigma$ \cite{Alikhanov2015A}
and weighted and shifted Gr\"{u}nwald difference (WSGD) \cite{tian2015class} formulae.
In Section 3, the B2T preconditioner and skew-circulant preconditioner are proposed and analyzed.
In Section 4, numerical examples are provided to demonstrate the efficiency of the two proposed preconditioners.
Some conclusions are drawn in Section 5.

\section{Finite difference discretization and the BLTT system}
\label{sec2}

In this section, the finite difference method is employed to discretize (\ref{eq1.1})
in both time and space. Then the BLTT system is derived
based on the obtained time-marching scheme.

\subsection{The time-marching scheme}
\label{sec2.1}

First of all, the WSGD operator is used to
approximate the left- and right- Riemann-Liouville derivatives \cite{tian2015class} (in this paper $(p,q) = (1,0)$).
Let $h = \frac{L}{N}$ be the grid spacing for the positive integer $N$. Hence the space domain is covered
by $\bar{\omega}_{h} = \{ x_i = i h | 0 \leq i \leq N \}$,
and approximations of the left- and right- Riemann-Liouville derivatives
	can be expressed respectively as:
\begin{equation}
	\sideset{_0}{^\beta_x}{\mathop{\mathcal{D}}} u(x,t)\mid_{x = x_i}
	\approx \frac{1}{h^{\beta}}\sum\limits_{k = 0}^{i +1}
	\omega_{k}^{(\beta)} u_{i - k + 1}, \qquad
	\sideset{_x}{^\beta_L}{\mathop{\mathcal{D}}} u(x,t)\mid_{x = x_i}
	\approx \frac{1}{h^{\beta}}\sum\limits_{k = 0}^{N - i +1}
	\omega_{k}^{(\beta)} u_{i + k - 1},
	\label{eq2.1}
\end{equation}
where $u_{i}$ is the numerical approximation to $u(x_i, t)$,
\begin{equation*}
	\omega_{0}^{(\beta)} = \frac{\beta}{2} g_{0}^{(\beta)}, \qquad
	\omega_{k}^{(\beta)} = \frac{\beta}{2} g_{k}^{(\beta)}
	+ \frac{2 - \beta}{2} g_{k - 1}^{(\beta)}, ~ k \geq 1
\end{equation*}
and
\begin{equation*}
	g_{0}^{(\beta)} = 1, \qquad
	g_{k}^{(\beta)} = \left( 1 - \frac{\beta + 1}{k} \right)
	g_{k - 1}^{(\beta)}, ~ k = 1, 2, \cdots.
\end{equation*}

Substituting Eq. (\ref{eq2.1}) into Eq. (\ref{eq1.1}),
the semi-discretized system of fractional ordinary differential equations is expressed as:
\begin{equation}
	\begin{cases}
		h^\beta \sideset{_0^C}{^\alpha_t}{\mathop{\mathcal{D}}} \bm{u}(t)
		= K_{N}\bm{u}(t) + h^\beta \bm{f}(t), & 0 < t \leq T, \\
		u(x,0) = u_{0}(x), & 0 \leq x \leq L,\\
	\end{cases}
	\label{eq2.2}
\end{equation}
where $\bm{u}(t) = \left[ u_{1}, u_{2}, \cdots, u_{N - 1} \right]^{T}$,
$\sideset{_0^C}{^\alpha_t}{\mathop{\mathcal{D}}} \bm{u}(t)
	= \left[ \sideset{_0^C}{^\alpha_t}{\mathop{\mathcal{D}}} u_1, \cdots,
	\sideset{_0^C}{^\alpha_t}{\mathop{\mathcal{D}}} u_{N - 1} \right]^{T}$,
$\bm{f}(t) = \left[ f_{1}, f_{2}, \cdots, f_{N - 1} \right]^{T}$ with $f_i = f(x_i,t)$ ($0 \leq i \leq N$) ,
$K_{N} = e_1 G_\beta + e_2 G_\beta^T$,
and the Toeplitz matrix $G_\beta$ is given
\begin{equation*}
	G_{\beta} =
	\begin{bmatrix}
		\omega_1^{(\beta)} &  \omega_0^{(\beta)} & 0 & \cdots & 0 & 0 \\
		\omega_2^{(\beta)} & \omega_1^{(\beta)} &  \omega_0^{(\beta)} & 0 & \cdots & 0 \\
		\vdots & \omega_2^{(\beta)} & \omega_1^{(\beta)} & \ddots & \ddots & \vdots \\
		\vdots & \ddots & \ddots & \ddots & \ddots & 0 \\
		\omega_{N-2}^{(\beta)} & \ddots & \ddots & \ddots & \omega_1^{(\beta)} & \omega_0^{(\beta)} \\
		\omega_{N-1}^{(\beta)} & \omega_{N-2}^{(\beta)} & \cdots & \cdots & \omega_2^{(\beta)} & \omega_1^{(\beta)}
	\end{bmatrix} \in \mathbb{R}^{(N - 1) \times (N - 1)}.
\end{equation*}

For a positive integer $M$,
the temporal partition is defined as
$\bar{\omega}_{\tau} = \{ t_j = j \tau, ~j = 0, 1, \cdots, M; ~t_M = T \}$
and let $u_i^j \approx u(x_i,t_j)$ be the approximate solution.
Through utilizing the $L2$-$1_\sigma$ formula \cite{Alikhanov2015A},
the temporal fractional derivative
$\sideset{_0^C}{^\alpha_t}{\mathop{\mathcal{D}}} u(x,t)$ can be discretized as:
\begin{equation}
	\sideset{_0^C}{^\alpha_t}{\mathop{\mathcal{D}}} u(x,t)\mid_{(x, t) = (x_i, t_{j + \sigma})}
	= \sum\limits^{j}_{s = 0} c^{(\alpha,\sigma)}_{j - s}
	\left( u_i^{s + 1} - u_i^{s} \right) + \mathcal{O}(\tau^{3 - \alpha}),
	\label{eq2.3}
\end{equation}
in which $\sigma = 1 - \alpha/2$ and for $j = 0$,
$c^{(\alpha,\sigma)}_{0} = \frac{\tau^{-\alpha}}{\Gamma(2 - \alpha)}  a^{(\alpha,\sigma)}_{0}$,
for $j \geq 1$,
\begin{equation*}
	c^{(\alpha,\sigma)}_{s}= \frac{\tau^{-\alpha}}{\Gamma(2 - \alpha)} \cdot
	\begin{cases}
		a^{(\alpha,\sigma)}_{0} + b^{(\alpha,\sigma)}_{1}, & s = 0,\\
		a^{(\alpha,\sigma)}_{s} + b^{(\alpha,\sigma)}_{s + 1} - b^{(\alpha,\sigma)}_{s}, & 1\leq s \leq j - 1,\\
		a^{(\alpha,\sigma)}_{j} - b^{(\alpha,\sigma)}_{j}, & s = j
	\end{cases}
\end{equation*}
with
\begin{equation*}
	a^{(\alpha,\sigma)}_{0} = \sigma^{1 - \alpha}, \qquad
	a^{(\alpha,\sigma)}_{l}
	= (l + \sigma)^{1 - \alpha} - (l - 1 + \sigma)^{1 - \alpha}~ ( l \geq 1),
\end{equation*}
\begin{equation*}
	b^{(\alpha,\sigma)}_{l} = \frac{1}{2 - \alpha}
	\left[(l + \sigma)^{2 - \alpha} - (l - 1 + \sigma)^{2 - \alpha} \right]
	- \frac{1}{2} \left[ (l + \sigma)^{1 - \alpha} - (l - 1 + \sigma)^{1 - \alpha} \right]~ ( l \geq 1).
\end{equation*}
Readers are suggested to refer to \cite{Alikhanov2015A} for a thoroughly discuss.

Substituting Eq. (\ref{eq2.3}) into Eq. (\ref{eq2.2}) and omitting the small term,
the discretized time-marching scheme is established as below
\begin{equation}
	h^\beta \sum\limits^{j}_{s = 0} c^{(\alpha,\sigma)}_{j - s}
	\left( \bm{u}^{s + 1} - \bm{u}^{s} \right)
	= K_N \bm{u}^{j + \sigma} + h^\beta \bm{f}^{j + \sigma},~j = 0, 1, \cdots, M - 1
	\label{eq2.4}
\end{equation}
with initial condition $u_i^0 = u_0(x_i)~(0 \leq i \leq N)$,
where $\bm{u}^{j + \sigma} = \sigma \bm{u}^{j + 1} + (1 - \sigma) \bm{u}^j$,
$\bm{u}^j = \left[ u_1^j, u_2^j, \cdots, u_{N - 1}^j \right]^T$,
$\bm{f}^{j + \sigma} = \left[ f_1^{j + \sigma}, f_2^{j + \sigma}, \cdots, f_{N - 1}^{j + \sigma} \right]^T$
and $f_i^{j + \sigma} = f(x_i, t_{j + \sigma})~(0 \leq i \leq N)$.
Furthermore, the stability and convergence with the second-order accuracy
of the time-marching scheme (\ref{eq2.4}) have been discussed in \cite{zhao2017fast}.

\subsection{The block lower triangular Toeplitz system}
\label{sec2.2}

Before deriving the BLTT system,
several auxiliary symbols are introduced:
$\bm{0}$ and $I$ represent zero and identity matrices of suitable orders, respectively.
$A_0 = h^\beta c_0^{(\alpha,\sigma)}I  - \sigma K_N$, $\bm{y}_0 = B \bm{u^0} + h^\beta \bm{f}^{\sigma}$,
\begin{equation*}
	A = \frac{\tau^{-\alpha} h^\beta}{\Gamma(2 - \alpha)} a_0^{(\alpha,\sigma)} I - \sigma K_N, \quad
	B = \frac{\tau^{-\alpha} h^\beta}{\Gamma(2 - \alpha)} a_0^{(\alpha,\sigma)} I + (1 - \sigma) K_N,
\end{equation*}
\begin{equation*}
	A_1 = h^\beta \left( c_1^{(\alpha,\sigma)} - c_0^{(\alpha,\sigma)} \right) I  - (1 - \sigma) K_N, \quad
	A_k = h^\beta \left( c_k^{(\alpha,\sigma)} - c_{k - 1}^{(\alpha,\sigma)} \right) I~(2 \leq k \leq M - 2).
\end{equation*}
To avoid misunderstanding, let
$v^{(\alpha,\sigma)}_{j} = \frac{\tau^{-\alpha}}{\Gamma(2 - \alpha)}
\left( a^{(\alpha,\sigma)}_{j} - b^{(\alpha,\sigma)}_{j} \right)$.
Then some other notations are given:
\begin{equation*}
	\bm{y}_1 = -\left[ h^\beta \left( v^{(\alpha,\sigma)}_{1} - c^{(\alpha,\sigma)}_{0} \right) I
	- (1 - \sigma) K_N \right] \bm{u}^1 + h^\beta \left( v^{(\alpha,\sigma)}_{1} \bm{u}^0
	+ \bm{f}^{1 + \sigma} \right),
\end{equation*}
\begin{equation*}
	\bm{y}_k = -h^\beta \left( v^{(\alpha,\sigma)}_{k} - c^{(\alpha,\sigma)}_{k - 1} \right) \bm{u}^1
	+ h^\beta \left( v^{(\alpha,\sigma)}_{k} \bm{u}^0 + \bm{f}^{k + \sigma} \right)~(2 \leq k \leq M - 1).
\end{equation*}

With the help of Eq. (\ref{eq2.4}), the BLTT system can be written as:
\begin{subnumcases}{}
	A \bm{u}^1 = \bm{y}_0,  \label{eq2.5a}\\
	W \bm{u} = \bm{y}, \label{eq2.5b}
\end{subnumcases}
where $\bm{y} = \left[ \bm{y}_1, \bm{y}_2, \cdots, \bm{y}_{M - 1} \right]^T$,
\begin{equation*}
	\bm{u} =
	\begin{bmatrix}
		\bm{u}^2 \\
		\bm{u}^3\\
		\vdots \\
		\bm{u}^{M}
	\end{bmatrix}, \quad
	W =
	\begin{bmatrix}
		A_0 & \bm{0} & \bm{0} & \cdots & \bm{0} \\
		A_1 & A_0 &  \bm{0} & \ddots & \vdots \\
		\vdots & \ddots & \ddots &  \ddots & \vdots \\
		A_{M - 3} & \ddots & \ddots & \ddots & \bm{0} \\
		A_{M - 2} & A_{M - 3} & \cdots &  \cdots & A_0
	\end{bmatrix}.
\end{equation*}

If the Kronecker product ``$\otimes$" is introduced, then Eq. (2.5) is equivalent to
\begin{equation*}
	\begin{cases}
		& A \bm{u}^1 = \bm{y}_0, \\
		& \tilde{W} \bm{u} = \bm{y},
	\end{cases}
\end{equation*}
in which $\tilde{W} = h^\beta \left( \tilde{A} \otimes I  \right)- \tilde{B} \otimes K_N$ with
\begin{equation*}
	\tilde{A} =
	\begin{bmatrix}
		c^{(\alpha,\sigma)}_{0} & 0 & 0 & \cdots & 0 & 0 \\
		c_1^{(\alpha,\sigma)} - c_0^{(\alpha,\sigma)} & c^{(\alpha,\sigma)}_{0} &  0 & 0 & \cdots & 0 \\
		\vdots & c_1^{(\alpha,\sigma)} - c_0^{(\alpha,\sigma)} & c^{(\alpha,\sigma)}_{0} & \ddots & \ddots & \vdots \\
		\vdots & \ddots & \ddots & \ddots & \ddots & 0 \\
		c_{M - 3}^{(\alpha,\sigma)} - c_{M - 4}^{(\alpha,\sigma)} & \ddots & \ddots & \ddots & c^{(\alpha,\sigma)}_{0} & 0 \\
		c_{M - 2}^{(\alpha,\sigma)} - c_{M - 3}^{(\alpha,\sigma)} & c_{M - 3}^{(\alpha,\sigma)} - c_{M - 4}^{(\alpha,\sigma)} &
		\cdots & \cdots &  c_1^{(\alpha,\sigma)} - c_0^{(\alpha,\sigma)} & c^{(\alpha,\sigma)}_{0}
	\end{bmatrix}
\end{equation*}
and $B = \textrm{tridiag}(1 - \sigma, \sigma, 0)$.

If the Gaussian elimination is adopted
for the BFS method \cite{huang2017fast} to solve (2.5),
the matrices $K_N$, $A$, $A_0$, $A_1$ and $B$ must be stored inherently.
Hence, the computational complexity and storage requirement of such the method are
	$\mathcal{O}(M N^3 + M N^2)$ and $\mathcal{O}(N^2)$, respectively.
To optimize the computational complexity,
we prefer to employ the preconditioned Krylov subspace methods to solve (2.5).
The key point of such preconditioned methods is to hunt for an efficient preconditioner.
In the following section, two economical preconditioners
are developed based on the special structures of $W$ and $A_0$,
and several properties of them are investigated.

\section{Two preconditioners and their spectra analysis}
\label{sec3}

In this section, two economical preconditioners are designed for solving Eq. (2.5).
The spectra of the preconditioned matrices are also analyzed.

\subsection{A block bi-diagonal Toeplitz preconditioner}
\label{sec3.1}
\begin{figure}[ht]
	\centering
	\subfigure[]
	{\includegraphics[width=3.0in,height=2.2in]{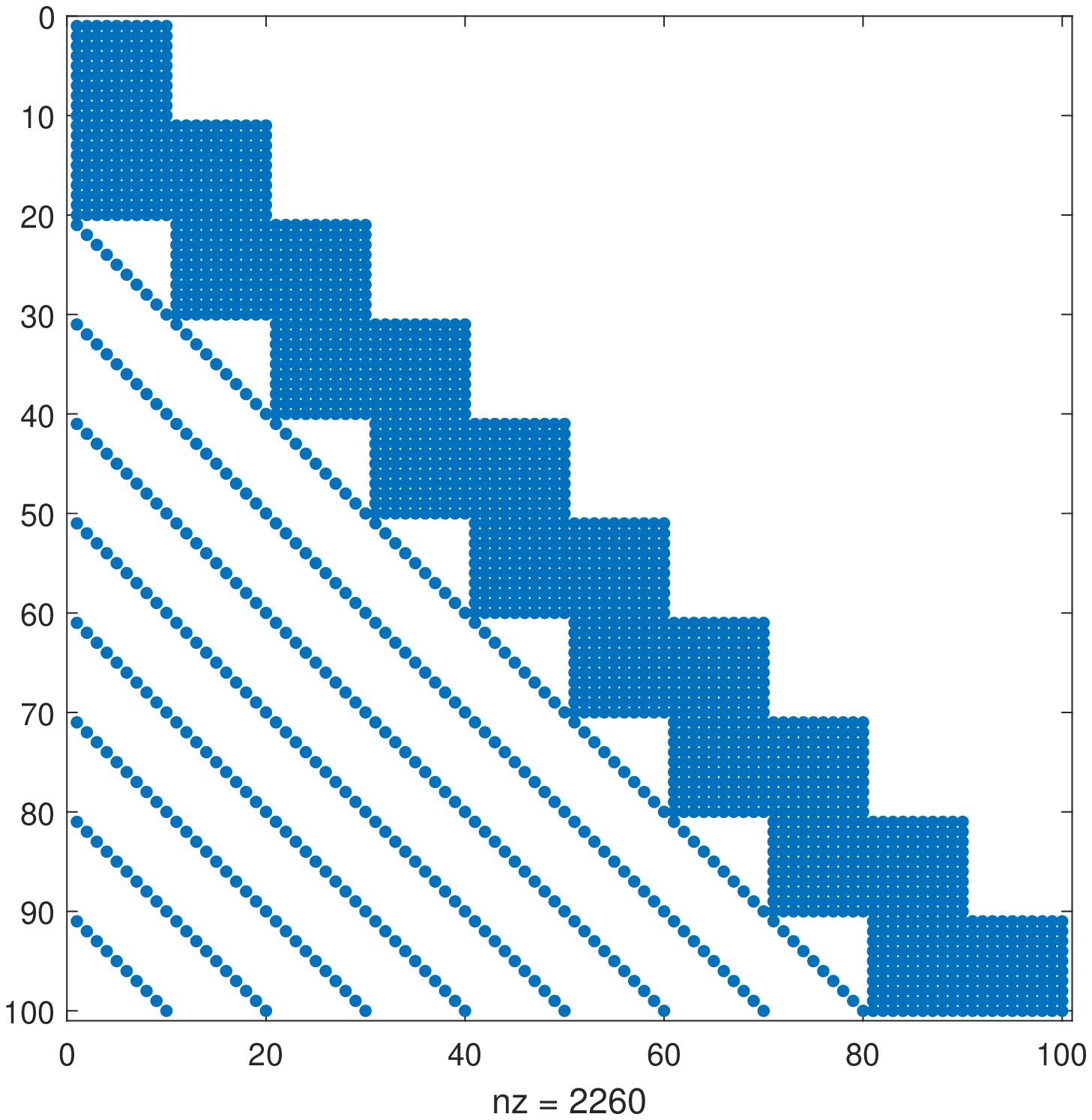}}
	\subfigure[]
	{\includegraphics[width=3.0in,height=2.2in]{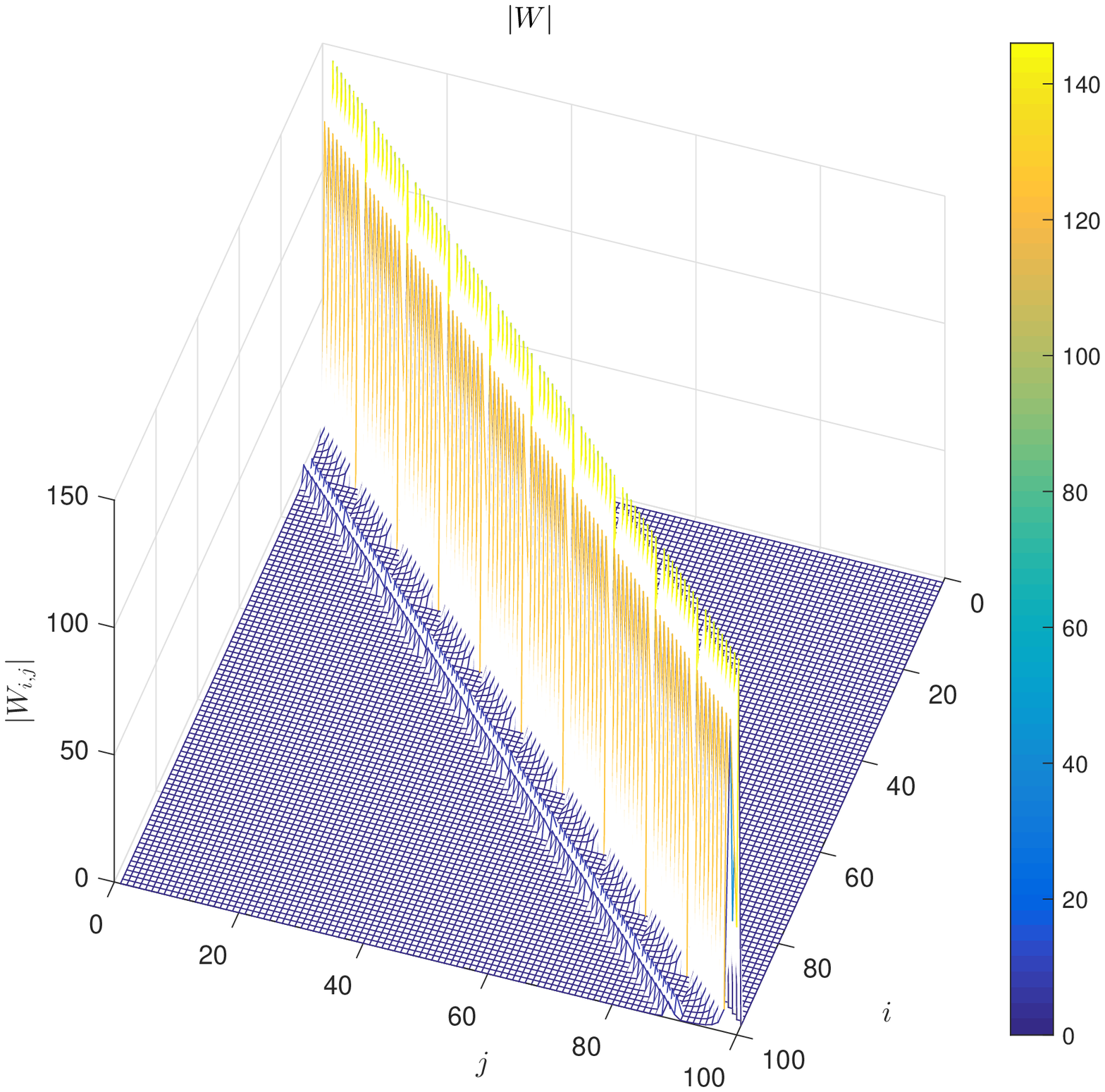}}
	\caption{The sparsity pattern (Left) and decay elements (Right)
		of matrix $W \in \mathbb{R}^{100 \times 100}$, when $M = N = 11$.}
	\label{fig1}
\end{figure}

To approximate the coefficient matrix $W$ well,
an example of the matrix $W$ is plotted in Fig. \ref{fig1} corresponding to $h = \tau = \frac{1}{11}$.
Fig. \ref{fig1}(a) shows the sparsity pattern of $W$.
From Fig. \ref{fig1}(b), it is noticeable that the diagonal entries of $W$ decay quickly,
i.e., the main information of $W$ clustered in the first two nonzero block diagonals.
Inspired by this observation, a block bi-diagonal Toeplitz preconditioner $P_W$
is developed for the linear system (\ref{eq2.5b}),
which only preserves the first two nonzero block diagonals of $W$, more precisely,
\begin{equation}
	P_W =
	\begin{bmatrix}
		A_0 & \bm{0} & \bm{0} & \cdots & \bm{0} & \bm{0} \\
		A_1 & A_0 & \bm{0} & \bm{0} & \cdots & \bm{0} \\
		\bm{0} & A_1 & A_0 & \ddots & \ddots & \vdots \\
		\vdots & \ddots & \ddots & \ddots & \ddots & \bm{0} \\
		\bm{0} & \ddots & \ddots & \ddots & A_0 & \bm{0} \\
		\bm{0} & \bm{0} & \cdots & \cdots & A_1 & A_0
	\end{bmatrix}.
	\label{eq3.1}
\end{equation}
Clearly $P_W$ is a block-Toeplitz matrix with Toeplitz-blocks,
thus its memory requirement is of $\mathcal{O}(N)$.

Several properties of $\omega^{(\beta)}_k$ are reviewed in the following lemma,
which is helpful to analyze the nonsingularity of $P_W$.
\begin{lemma}(\cite{zhao2017fast, feng2016high})
	Suppose that $1< \beta <2$, then the coefficients $\omega_{k}^{(\beta)}$ satisfy
	\begin{equation*}
		\begin{cases}
			\omega_{0}^{(\beta)} = \frac{\beta}{2} > 0, ~~
			\omega_{1}^{(\beta)} = \frac{2 - \beta - \beta^2}{2} < 0, ~~
			\omega_{2}^{(\beta)} = \frac{\beta (\beta^2 + \beta - 4)}{4}, \\
			1 \geq \omega_{0}^{(\beta)} \geq \omega_{3}^{(\beta)} \geq \omega_{4}^{(\beta)} \geq \cdots \geq 0, ~~
			\omega_{0}^{(\beta)} + \omega_{2}^{(\beta)} > 0, \\
			\sum\limits_{k = 0}^{\infty} \omega_{k}^{(\beta)} = 0, ~~
			\sum\limits_{k = 0}^{N} \omega_{k}^{(\beta)} < 0, ~~N \geq 2.
		\end{cases}
	\end{equation*}
	\label{lemma3.1}
\end{lemma}
As seen from Lemma \ref{lemma3.1}, we proceed to analyze the nonsingularity of $P_W$.
\begin{theorem}
	$P_W$ given in (\ref{eq3.1}) is nonsingular.
	\label{th3.1}
\end{theorem}
\textbf{Proof.}
Since $P_W$ is a block lower bi-diagonal matrix, the proof of  this theorem
is equivalent to prove the nonsingularity of $A_0$.

Firstly, we show that all eigenvalues of matrix $H = \frac{K_N + K_N^T}{2}$ are strictly negative.
From the definition of $K_N$ in (\ref{eq2.2}),
it has $H = \frac{e_1 + e_2}{2} \left( G_\beta + G_\beta^T \right)$.
Then according to the Gershgorin circle theorem \cite{varga2004gersgorin},
the $i$-th Gershgorin disc of $H $ is centered at
$\left( e_1 + e_2 \right) \omega_1^{(\beta)} < 0$ with radius
\begin{equation*}
	r_i^\beta = \frac{e_1 + e_2}{2} \left( \sum\limits_{k = 0, k \neq 1}^{i} \omega_k^{(\beta)}
	+  \sum\limits_{k = 0, k \neq 1}^{N - i} \omega_k^{(\beta)} \right)
	\leq \left( e_1 + e_2 \right) \sum\limits_{k = 0, k \neq 1}^{N} \omega_k^{(\beta)}
	< -\left( e_1 + e_2 \right) \omega_1^{(\beta)} ~ (1 \leq i \leq N - 1),
\end{equation*}
in which Lemma \ref{lemma3.1} is adopted.
Thus the real parts of all eigenvalues of $A_0$ are strictly positive.
The proof of Theorem \ref{th3.1} is completed.
\hfill $\Box$

Theorem \ref{th3.1} also implies that the matrices $A$ and $W$ are invertible.
Now, the eigenvalues of the preconditioned matrix $P_W^{-1} W$ can be studied.
\begin{theorem}
	The eigenvalues of the preconditioned matrix $P_W^{-1} W$ are all equal to $1$.
	\label{th3.2}
\end{theorem}
\textbf{Proof.}
It is known that the product of two block lower triangular matrices
also is a block lower triangular matrix.
After simple calculations, it notes that
\begin{equation*}
	P_W^{-1} W =
	\begin{bmatrix}
		I & \bm{0} & \cdots & \cdots & \bm{0} \\
		\bm{0} & I &  \ddots & \ddots & \vdots \\
		J_2 & \ddots & \ddots &  \ddots & \vdots \\
		\vdots & \ddots & \ddots & \ddots & \bm{0} \\
		J_{M - 2} & J_{M - 3} & \cdots &  \cdots & I
	\end{bmatrix}
\end{equation*}
is a block lower triangular matrix,
where $J_2 = A_0^{-1} A_2$, $J_k = A_0^{-1} \left( A_k - A_1 J_{k - 1} \right)~(3 \leq k \leq M - 2)$.
From the above equality, the main diagonal elements of $P_W^{-1} W$ are $1$,
which completes the proof. \hfill $\Box$
\begin{remark}
	The preconditioned Krylov subspace methods require us to compute $P_W^{-1} \bm{v}$, where $\bm{v}$ is a vector.
	In this work, the Thomas method is employed to compute such matrix-vector multiplications.
	Hence, only $A_0^{-1} \bm{v}$ is needed to compute.
	In practical computation, the Toeplitz inversion formula \cite{ng2003recursive}
	combined with Krylov subspace methods is used to calculate $A_0^{-1} \bm{v}$,
	and this will be discussed in Section \ref{sec3.2}.
	\label{remark1}
\end{remark}

For the sake of clarity, the Thomas method for calculating $P_W^{-1} \bm{v}$ is given as below.
\begin{algorithm}[H]
	\caption{Compute $\bm{z} = P_W^{-1} \bm{v}$}
	\begin{algorithmic}[1]
		\STATE Reshape $\bm{v}$ into an $(N - 1) \times M$ matrix $\check{V}$
		\STATE Compute $\bm{\hat{b}}_1 = A_0^{-1} \check{V}(:,1)$ via Algorithm \ref{alg2} in Section \ref{sec3.2}
		\FOR {$k = 2, \cdots, M$}
		\STATE $\bm{\varphi} = \check{V}(:,k) - A_1 \hat{\bm{b}}_{k - 1}$
		\STATE $\bm{\hat{b}}_k = A_0^{-1} \bm{\varphi}$ via Algorithm \ref{alg2} in Section \ref{sec3.2}
		\ENDFOR
		\STATE Stack $\bm{\hat{b}}_k~(k = 1,\cdots,M)$ in a vector $\bm{z}$
	\end{algorithmic}
	\label{alg1}
\end{algorithm}

In line 4 of Algorithm \ref{alg1}, the matrix-vector multiplications can be done
via FFTs in $\mathcal{O}(N \log N)$ operations \cite{ng2004iterative,chan2007introduction}.
As for the storage requirement, $\bm{v}$, $\bm{\hat{b}}_k$, $\bm{\varphi}$,
the first column and first row of $A_1$ must be stored.
Thus only $\mathcal{O}(M N)$ memory is needed in Algorithm \ref{alg1}.

\subsection{A skew-circulant preconditioner}
\label{sec3.2}

According to the Toeplitz inversion formula in \cite{ng2003recursive}, two Toeplitz systems
\begin{equation}
	\begin{cases}
		A_0 \bm{\xi} = \bm{q}_1, \\
		A_0 \bm{\eta} = \bm{q}_{N - 1}
	\end{cases}
	\label{eq3.2}
\end{equation}
require to be solved, where $\bm{\xi} = [\xi_1, \cdots,\xi_{N - 1}]^T$,
$\bm{\eta} = [\eta_1, \cdots,\eta_{N - 1}]^T$, $\bm{q}_1$ and $\bm{q}_{N - 1}$
are the first and last columns of the identity matrix of order $(N - 1)$, respectively.
As mentioned in Remark \ref{remark1}, Krylov subspace methods are chosen to solve (\ref{eq3.2}).
However, when $A_0$ is ill-conditioned, Krylov subspace methods converge very slowly.
To remedy such difficulties, in this subsection,
a new skew-circulant preconditioner $P_{sk}$ is designed
and the spectrum of $P_{sk}^{-1} A_0$ is discussed.
The expression of our skew-circulant $P_{sk}$ is given as follows
\begin{equation}
	P_{sk} = h^\beta c_0^{(\alpha,\sigma)}I  - \sigma sk(K_N),
	\label{eq3.3}
\end{equation}
where $sk(K_N) = e_1 sk(G_\beta) + e_2 sk(G_\beta)^T$ with
\begin{equation*}
	sk(G_{\beta}) =
	\begin{bmatrix}
		\omega_1^{(\beta)} &  \omega_0^{(\beta)} & -\omega_{N-2}^{(\beta)} & \cdots & -\omega_{2}^{(\beta)} \\
		\omega_2^{(\beta)} & \omega_1^{(\beta)} &  \omega_0^{(\beta)} & \ddots & \vdots \\
		\vdots & \ddots & \ddots &  \ddots & -\omega_{N-2}^{(\beta)} \\
		\omega_{N-2}^{(\beta)} & \ddots & \ddots & \ddots & \omega_0^{(\beta)} \\
		-\omega_{0}^{(\beta)} & \omega_{N-2}^{(\beta)} & \cdots &  \omega_2^{(\beta)} & \omega_1^{(\beta)}
	\end{bmatrix} \in \mathbb{R}^{(N - 1) \times (N - 1)}.
\end{equation*}
Similar to the proof of Theorem \ref{th3.1},
the following theorem provide an essential property of $P_{sk}$ in (\ref{eq3.3}).
\begin{theorem}
	The matrix $P_{sk}$ given in (\ref{eq3.3}) is invertible.
	\label{th3.3}
\end{theorem}
\textbf{Proof.}
Firstly, we prove that all eigenvalues of matrix $\hat{H} = -\frac{sk(K_N) + sk(K_N)^T}{2}$
are strictly positive. Based on the definition of $sk(G_\beta)$
and the Gershgorin circle theorem \cite{varga2004gersgorin},
all the Gershgorin disc of the matrix $\hat{H}$ are centered at
$-\left( e_1 + e_2 \right) \omega_1^{(\beta)} > 0$
with radius
\begin{equation*}
	r = \frac{e_1 + e_2}{2} \left[ 2\left( \omega_0^{(\beta)} + \omega_2^{(\beta)} \right)
	+ \sum\limits_{k = 3}^{N - 2} \left| \omega_k^{(\beta)} - \omega_{N + 1 - k}^{(\beta)} \right| \right]
	\leq \left( e_1 + e_2 \right) \sum\limits_{k = 0, k \neq 1}^{N} \omega_k^{(\beta)}
	\leq -\left( e_1 + e_2 \right) \omega_1^{(\beta)}.
\end{equation*}
Thus, the real parts of all eigenvalues of $P_{sk}$ are strictly positive.
Then the targeted result follows. \hfill $\Box$

An $n \times n$ skew-circulant matrix $\mathcal{C}$ has the spectral decomposition
\cite{ng2004iterative, chan2007introduction}:
\begin{equation*}
	\mathcal{C} = \Omega^* F^* \Lambda F \Omega,
\end{equation*}
here $\Omega = \textrm{diag} \left[ 1, (-1)^{-1/n}, \cdots, (-1)^{-(n - 1)/n} \right]$,
$F$ is the discrete Fourier matrix, $F^*$ represents the conjugate transpose of $F$,
and $\Lambda$ is a diagonal matrix containing all eigenvalues of $\mathcal{C}$.
Let $sk(G_\beta) = \Omega^* F^* \Lambda_s F \Omega$, then
$sk(G_\beta)^T = \Omega^* F^* \bar{\Lambda}_s F \Omega$ and
$P_{sk} = \Omega^* F^* \Lambda F \Omega$,
where $\Lambda = h^\beta c_0^{(\alpha,\sigma)}I  - \sigma \left( e_1 \Lambda_s + e_2 \bar{\Lambda}_s \right)$
and $\bar{\Lambda}_s$ is the complex conjugate of $\Lambda_s$.
With the help of the decomposition of $P_{sk}$, the following result is obtained immediately.
\begin{lemma}
	Suppose $0 < \hat{v} < h^\beta c_0^{(\alpha,\sigma)}$,
	then $\parallel P_{sk}^{-1} \parallel_2 \leq \frac{1}{\hat{v}}$.
	\label{lemma3.2}
\end{lemma}
\textbf{Proof.}
By Theorem \ref{th3.3}, we obtain $Re([\Lambda_s]_{k,k}) < 0$,
where $Re([\Lambda_s]_{k,k})$ means the real part of $[\Lambda_s]_{k,k}$.
Then
\begin{equation*}
	\left| [\Lambda]_{k,k} \right| \geq Re([\Lambda]_{k,k})
	= h^\beta c_0^{(\alpha,\sigma)}  - \sigma \left( e_1 Re([\Lambda_s]_{k,k})
	+ e_2 Re([\bar{\Lambda}_s]_{k,k}) \right)
	\geq \hat{v}, ~k = 1, 2, \cdots, N - 1.
\end{equation*}
Therefore
\begin{equation*}
	\parallel P_{sk}^{-1} \parallel_2
	= \frac{1}{\min\limits_{1\leq k\leq N - 1} \left| [\Lambda]_{k,k} \right| } \leq \frac{1}{\hat{v}}.
\end{equation*}
\hfill $\Box$

To analyze the spectrum of $P_{sk}^{-1} A_0$,
we first prove that the generating function of the Toeplitz matrix $K_N$ is
in the Wiener class \cite{chan2007introduction}.
\begin{lemma}
	The generating function of the sequence $\left\{ K_N \right\}_{N = 2}^{\infty}$ is
	in the Wiener class.
	\label{lemma3.3}
\end{lemma}
\textbf{Proof.}
For the Toeplitz matrix $K_N$ in (\ref{eq2.2}), its generating function is
\begin{equation*}
	p(\theta) = \sum\limits_{k = -\infty}^{\infty} \ell_k e^{\bm{i} k \theta}
	= \sum\limits_{k = -1}^{\infty} \omega_{k + 1}^{(\beta)}
	\left( e_1 e^{\bm{i} k \theta} + e_2 e^{-\bm{i} k \theta} \right),
\end{equation*}
where $\bm{i} = \sqrt{-1}$ and $\theta \in [-\pi, \pi]$.
By the properties of $\omega_k^{(\beta)}$, it yields
\begin{equation*}
	\sum\limits_{k = -\infty}^{\infty} \left| \ell_k \right|
	\leq \left( e_1 + e_2 \right) \sum\limits_{k = -1}^{\infty} \left| \omega_{k + 1}^{(\beta)}  \right|
	= \left( e_1 + e_2 \right)
	\left(- 2\omega_{1}^{(\beta)} + \left| \omega_{2}^{(\beta)} \right| - \omega_{2}^{(\beta)} \right)  < \infty.
\end{equation*}
Thus, the generating function $p(\theta)$ is in the Wiener class.
\hfill $\Box$

According to Lemma \ref{lemma3.3}, the following result is true.
\begin{lemma}
	Let $p(\theta)$ be the generating function of $K_N$.
	Then for any $\varepsilon > 0$, there exists an $N' > 0$, such that for all $N > N' + 1$,
	$A_0 - P_{sk} = \tilde{U} + \tilde{V}$, where $rank(\tilde{U}) < 2N'$
	and $\parallel \tilde{V} \parallel_2 \leq \varepsilon$.
	\label{lemma3.4}
\end{lemma}
\textbf{Proof.}
Define $D_{sk} = A_0 - P_{sk} = \sigma \left( sk(K_N) - K_N \right)$.
It can be checked that $D_{sk}$ is a Toeplitz matrix,
and its first column and first row are respectively
\begin{equation*}
	\begin{split}
		&-\sigma [0, 0, e_2 \omega_{N - 2}^{(\beta)}, \cdots, e_2 \omega_{3}^{(\beta)},
		e_1 (\omega_{0}^{(\beta)} + \omega_{N - 2}^{(\beta)}) + e_2 \omega_{2}^{(\beta)} ]^T,\\
		&-\sigma [0, 0, e_1\omega_{N - 2}^{(\beta)}, \cdots, e_1 \omega_{3}^{(\beta)},
		e_1 \omega_{2}^{(\beta)} + e_2 (\omega_{0}^{(\beta)} + \omega_{N - 2}^{(\beta)})].
	\end{split}
\end{equation*}

Using Lemma \ref{lemma3.3}, we know that $p(\theta)$ is in the Wiener class.
Then for any $\varepsilon > 0$, there exists an $N' > 0$
such that $\sum\limits_{k = N' + 1}^{\infty} \left| \ell_k \right|
= e_2 \sum\limits_{k = N' + 1}^{\infty} \left| \omega_{k + 1}^{(\beta)} \right|
\leq \frac{ e_2 \varepsilon}{\sigma (e_1 + e_2)}$.
Let $\tilde{V}$ be the $(N - 1)$-by-$(N - 1)$ matrix obtained from $D_{sk}$ by copying the
$(N - 1 - N')$-by-$(N - 1 - N')$ leading principal submatrix of $D_{sk}$.
Hence the leading $(N - 1 - N') \times (N - 1 - N')$ block of $\tilde{V}$ is a Toeplitz matrix.
Thus
\begin{equation*}
	\begin{split}
		\parallel \tilde{V} \parallel_1
		& = \sigma \max\left\{ \sum\limits_{k = N' + 1}^{N - 3} \left| \ell_{k} \right|,
		\sum\limits_{k = N' + 1}^{N - 3} \left| \ell_{-k} \right|,
		\max\limits_{3 \leq j \leq N - 3 - N'} \left( \sum\limits_{k = N' + j}^{N - 3} \left| \ell_k \right|
		+ \sum\limits_{k = N - j}^{N - 3} \left| \ell_{-k} \right| \right) \right\}\\
		& \leq \sigma (e_1 + e_2) \sum\limits_{k = N' + 1}^{\infty} \left| \omega_{k + 1}^{(\beta)} \right|
		\leq \varepsilon.
	\end{split}
\end{equation*}
Similarly, $\parallel \tilde{V} \parallel_\infty \leq \varepsilon$.
Thus $\parallel \tilde{V} \parallel_2 \leq
\left( \parallel \tilde{V} \parallel_1 \cdot \parallel \tilde{V} \parallel_\infty \right)^{1/2} \leq \varepsilon$.

Let $\tilde{U} = D_{sk} - \tilde{V}$. It is obvious that $\tilde{U}$ is an $(N - 1) \times (N - 1)$ matrix obtained
from $D_{sk}$ by replacing the $(N - 1 - N') \times (N - 1 - N')$ leading principal submatrix of $D_{sk}$ by the
zero matrix. Hence $rank(\tilde{U}) \leq 2 N'$.
\hfill $\Box$

Combining Lemmas \ref{lemma3.2} and \ref{lemma3.4}, the spectrum of $P_{sk}^{-1} A_0 - I$
is discussed.
\begin{theorem}
	Suppose $0 < \hat{v} < h^\beta c_0^{(\alpha,\sigma)}$.
	Then for any $\varepsilon > 0$, there exists an $N' > 0$, such that for all $N - 1 > N'$,
	$P_{sk}^{-1} A_0 - I = U + V$, where $rank(U) < 2N'$
	and $\parallel V \parallel_2 \leq \frac{\varepsilon}{\hat{v}}$.
	\label{th3.4}
\end{theorem}
\textbf{Proof.}
According to Lemma \ref{lemma3.4}, for any $\varepsilon > 0$, there exists an $N' > 0$,
such that for all $N - 1 > N'$,
\begin{equation*}
	P_{sk}^{-1} A_0 - I
	= P_{sk}^{-1} \left( A_0 - P_{sk} \right)
	= U + V,
\end{equation*}
where $U = P_{sk}^{-1} \tilde{U}$ and $V = P_{sk}^{-1} \tilde{V}$.
Applying Lemma \ref{lemma3.2}, it yields
\begin{equation*}
	\parallel V \parallel_2
	= \parallel P_{sk}^{-1} \tilde{V} \parallel_2
	\leq  \parallel P_{sk}^{-1} \parallel_2 \parallel \tilde{V} \parallel_2
	\leq \frac{\varepsilon}{\hat{v}}.
\end{equation*}
On the other hand, $rank(U) =  rank(P_{sk}^{-1} \tilde{U}) < 2N'$.
\hfill $\Box$

\begin{remark}
	Since the matrix $A$ is slightly different to $A_0$,
	the $P_{sk}$ in (\ref{eq3.3}) still works for solving (\ref{eq2.5a}).
	Hence, in this work, $P_{sk}$ is also applied to solve (\ref{eq2.5a}).
\end{remark}

For convenience, our strategy in this subsection is concluded in the following algorithm.
\begin{algorithm}
	\caption{Compute $\bm{\tilde{z}} = A_0^{-1} \bm{v}$}
	\begin{algorithmic}[1]
		\STATE {Solve $A_0 \bm{\xi} = \bm{q}_1$ via FGMRES/PBiCGSTAB with $P_{sk}$ \\
			Solve $A_0 \bm{\eta} = \bm{q}_{N - 1}$ via FGMRES/PBiCGSTAB with $P_{sk}$}
		\STATE {$\bm{s}_1 = [\eta_{N - 1}, -\eta_1,\cdots, -\eta_{N - 2}]^T$,
			$\bm{s}_2 = [\eta_{N - 1}, \eta_1,\cdots, \eta_{N - 2}]^T$}
		\STATE {$\Lambda^{(1)} = \textrm{diag}(F \bm{\xi})$, ~$\Lambda^{(2)} = \textrm{diag}(\Omega^{*} F \bm{s}_1)$, \\
			$\Lambda^{(3)} = \textrm{diag}(F \bm{s}_2)$, $\Lambda^{(4)} = \textrm{diag}(\Omega^{*} F \bm{\xi})$}
		\STATE {$\bm{\tilde{v}} = F \Omega \bm{v}$}
		\STATE {$\bm{z}_1 = \Omega^{*} F^{*} \Lambda^{(2)}\bm{\tilde{v}}$,
			~$\bm{z}_2 = \Omega^{*} F^{*} \Lambda^{(4)}\bm{\tilde{v}}$, \\
			$\bm{z}_3 = \Lambda^{(1)} F \bm{z}_1$, \quad~ $\bm{z}_4 = \Lambda^{(3)} F \bm{z}_2$}
		\STATE $\bm{\tilde{z}} = \frac{1}{2 \xi_1} F^{*} (\bm{z}_3 + \bm{z}_4)$
	\end{algorithmic}
	\label{alg2}
\end{algorithm}

In this algorithm, ten fast Fourier transforms are needed.
Thus, the complexity and storage requirement
are $\mathcal{O}(N \log N)$ and $\mathcal{O}(N)$, respectively.

\section{Numerical experiments}
\label{sec4}

In this section, one example is reported to show the performance
of the proposed preconditioners in Section \ref{sec3}.
In order to illustrate the efficiency of $P_{sk}$,
the Strang's circulant preconditioner \cite{ng2004iterative,chan2007introduction} is also tested,
which can be written as
\begin{equation*}
	P_{s} = h^\beta c_0^{(\alpha,\sigma)}I  - \sigma s(K_N),
\end{equation*}
where $s(K_N) = e_1 s(G_\beta) + e_2 s(G_\beta)^T$.
More precisely, the first columns of circulant matrices $s(G_\beta)$ and $s(G_\beta)^T$ are
$\left[ \omega_1^{(\beta)}, \cdots, \omega_{\lfloor N/2 \rfloor}^{(\beta)}, 0,
\cdots, 0, \omega_0^{(\beta)} \right]^T$ and
$\left[ \omega_1^{(\beta)}, \omega_0^{(\beta)}, 0, \cdots, 0,
\omega_{\lfloor N/2 \rfloor}^{(\beta)}, \cdots, \omega_2^{(\beta)} \right]^T$,  respectively.

The PBiCGSTAB and FGMRES methods for solving (2.5) terminate
if the relative residual error satisfies
$\frac{\| \bm{r}^{(k)} \|_2}{\| \bm{r}^{(0)} \|_2} < 10^{-8}$
or the iteration number is more than $1000$, where $\bm{r}^{(k)}$ denotes residual vector in the $k$-th iteration,
and the initial guess is chosen as the zero vector. Since the $P_W$ as a preconditioner for solving (2.5),
it is not necessary to compute the $P^{-1}_W \bm{v}$ accurately.
Hence the stop criterion of PBiCGSTAB or FGMRES methods in Algorithm \ref{alg2} is
$\frac{\| \bm{r}^{(k)} \|_2}{\| \bm{r}^{(0)} \|_2} < 10^{-3}$,
and the initial guess is also chosen as the zero vector. All of the symbols shown below will appear in later.

All experiments are carried out via MATLAB 2017a on a Windows 10 (64 bit) PC with the configuration:
Intel(R) Core(TM) i7-7700T CPU 2.90 GHz and 8 GB RAM.
\begin{center}\footnotesize\tabcolsep=2.5pt
	\begin{tabular}{|l|l|}
		\hline
		Symbol & Explanation \\
		\hline
		BS &  The MATLAB's backslash method to solve (2.5) \\
		BFSM & The BFS method to solve (2.5) \\
		SK2-PBiCGSTAB & The PBiCGSTAB method with the preconditioners $P_W$ and $P_{sk}$ to solve (2.5) \\
		SK2-FGMRES & The FGMRES method with the preconditioners $P_W$ and $P_{sk}$ to solve (2.5) \\
		S2-PBiCGSTAB & The PBiCGSTAB method with the preconditioners $P_W$ and $P_{s}$ to solve (2.5) \\
		S2-FGMRES & The FGMRES method with the preconditioners $P_W$ and $P_{s}$ to solve (2.5) \\
		$\mathrm{Iter1}$ & The number of iterations required for solving (\ref{eq2.5a}) \\
		$\mathrm{Iter2}$ & The number of iterations required for solving (\ref{eq2.5b}) \\
		$\mathrm{Iter3}$ & The number of iterations required for solving (\ref{eq3.2}) \\
		$\mathrm{Iter}$ & $\mathrm{Iter1} + \mathrm{Iter2}$ \\
		$\textrm{Time}$ & Total CPU time in seconds for solving the whole BLTT system (2.5) \\
		\dag & Out of memory \\
		\hline
	\end{tabular}
\end{center}

\begin{figure}[ht]
	\centering
	\subfigure[Eigenvalues of $W$]
	{\includegraphics[width=3.0in,height=2.2in]{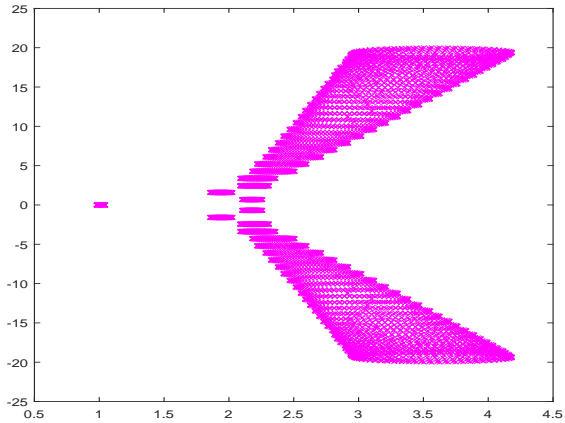}}
	\subfigure[Eigenvalues of $P_{W}^{-1} W$]
	{\includegraphics[width=3.0in,height=2.25in]{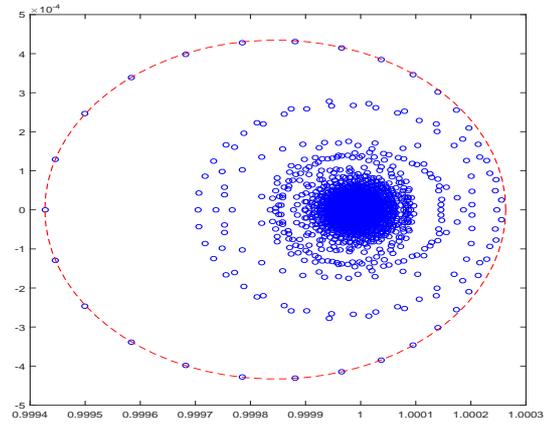}}\\
	\subfigure[Eigenvalues of $W$]
	{\includegraphics[width=3.0in,height=2.2in]{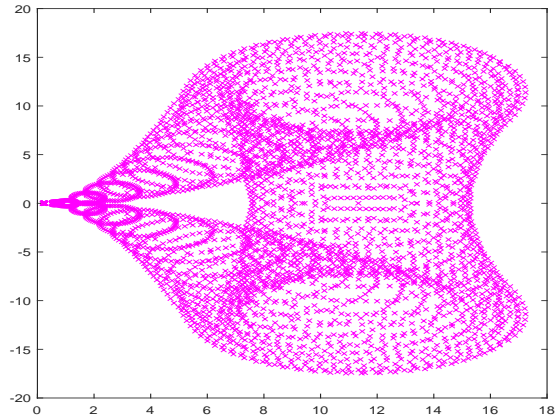}}
	\subfigure[Eigenvalues of $P_{W}^{-1} W$]
	{\includegraphics[width=3.0in,height=2.2in]{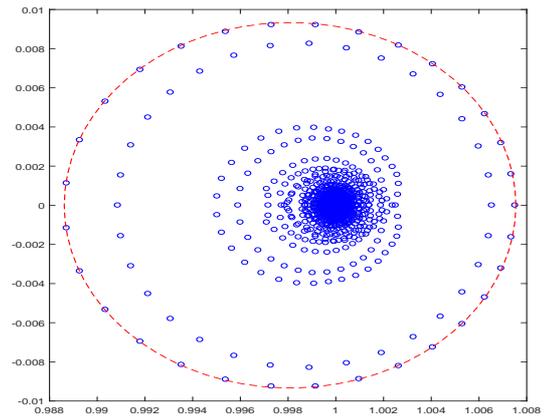}}
	\caption{Spectra of $W$ and $P_{W}^{-1} W$, when $M = N = 2^6$ in Example 1.
		Top row: $(\alpha, \beta) = (0.1, 1.1)$; Bottom row: $(\alpha, \beta) = (0.7, 1.4)$.}
	\label{fig2}
\end{figure}

\noindent \textbf{Example 1.}
Considering Eq. (\ref{eq1.1}) with diffusion coefficients $e_1 = 20$, $e_2 = 0.02$,
the source term
\begin{equation*}
	\begin{split}
		f(x,t) = & 2 t^{1 - \alpha} E_{1, 2 - \alpha}(2t) x^2 (1 - x)^2
		- e^{2t} \Bigg\{ \frac{\Gamma(3)}{\Gamma(3 - \beta)} \left[ e_1 x^{2 - \beta} + e_2 (1 - x)^{2 - \beta} \right] \\
		& - \frac{2 \Gamma(4)}{\Gamma(4 - \beta)} \left[ e_1 x^{3 - \beta} + e_2 (1 - x)^{3 - \beta} \right]
		+ \frac{\Gamma(5)}{\Gamma(5 - \beta)} \left[ e_1 x^{4 - \beta} + e_2 (1 - x)^{4 - \beta} \right] \Bigg\},
	\end{split}
\end{equation*}
in which $E_{\mu, \nu}(z)$ is the Mittag-Leffler function \cite{podlubny1998fractional} with two parameters defined by
\begin{equation*}
	E_{\mu, \nu}(z) = \sum_{k = 0}^{\infty} \frac{z^{k}}{\Gamma(\mu k + \nu)}.
\end{equation*}
The exact solution of the TSFDE problem (\ref{eq1.1}) is $u(x,t) = e^{2 t} x^2 (1 - x)^2$.
\begin{figure}[ht]
	\centering
	\subfigure[Eigenvalues of $A_0$]
	{\includegraphics[width=3.0in,height=2.2in]{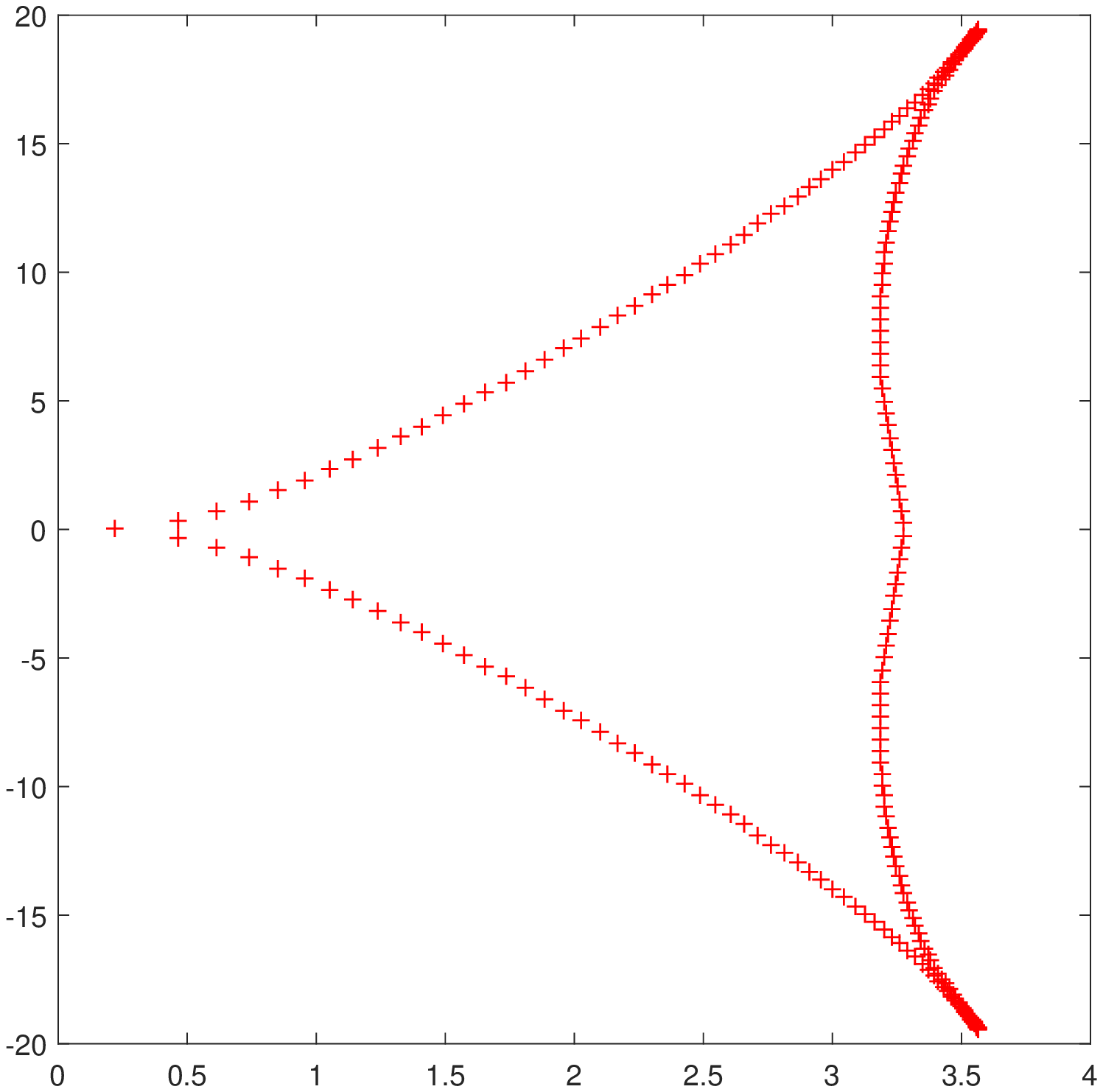}}
	\subfigure[Eigenvalues of $P_{sk}^{-1} A_0$ ({\color{blue} $*$})
	and $P_{s}^{-1} A_0$ ({\color{magenta} $\times$})]
	{\includegraphics[width=3.0in,height=2.2in]{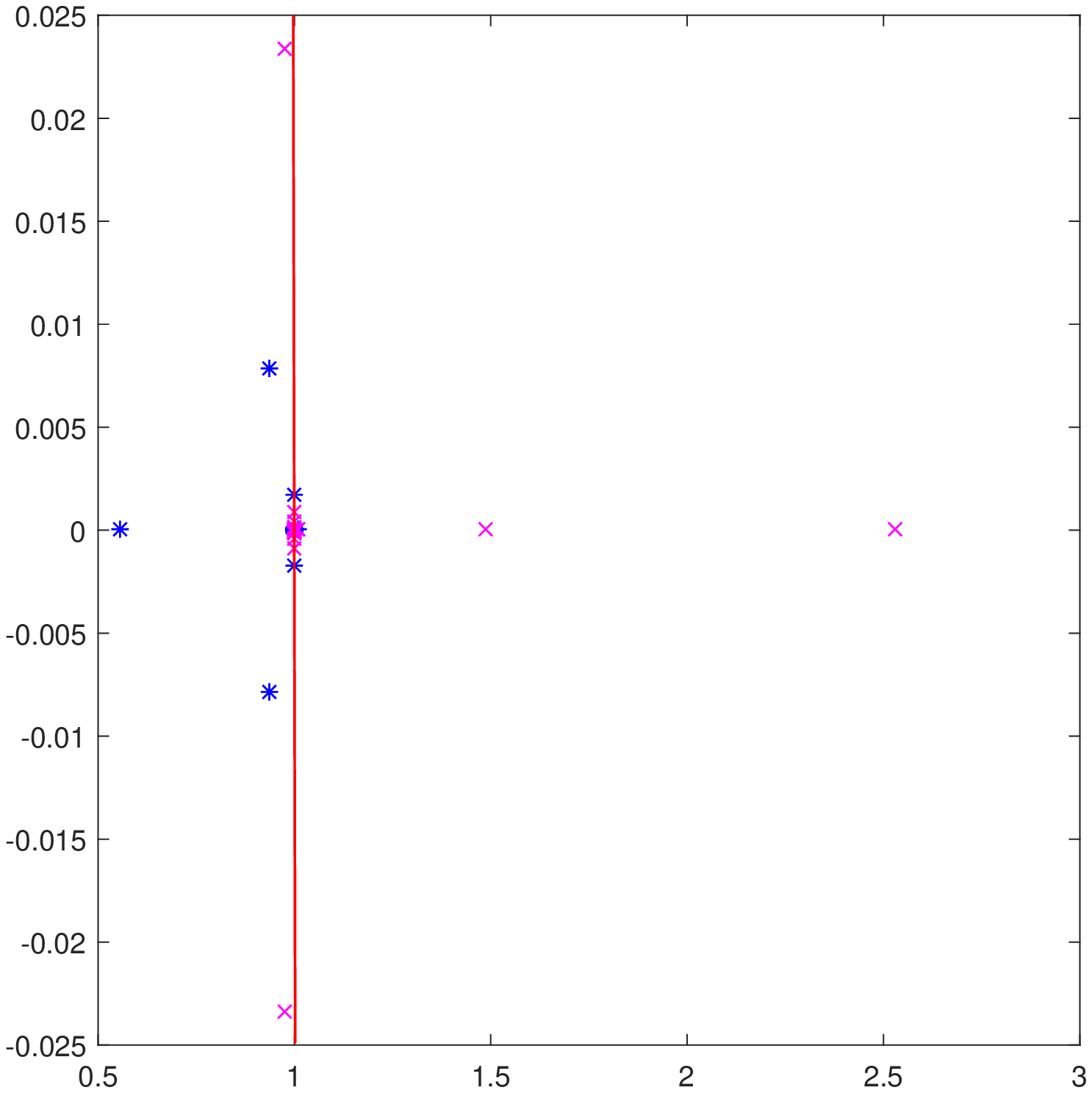}}\\
	\subfigure[Eigenvalues of $A_0$]
	{\includegraphics[width=3.0in,height=2.2in]{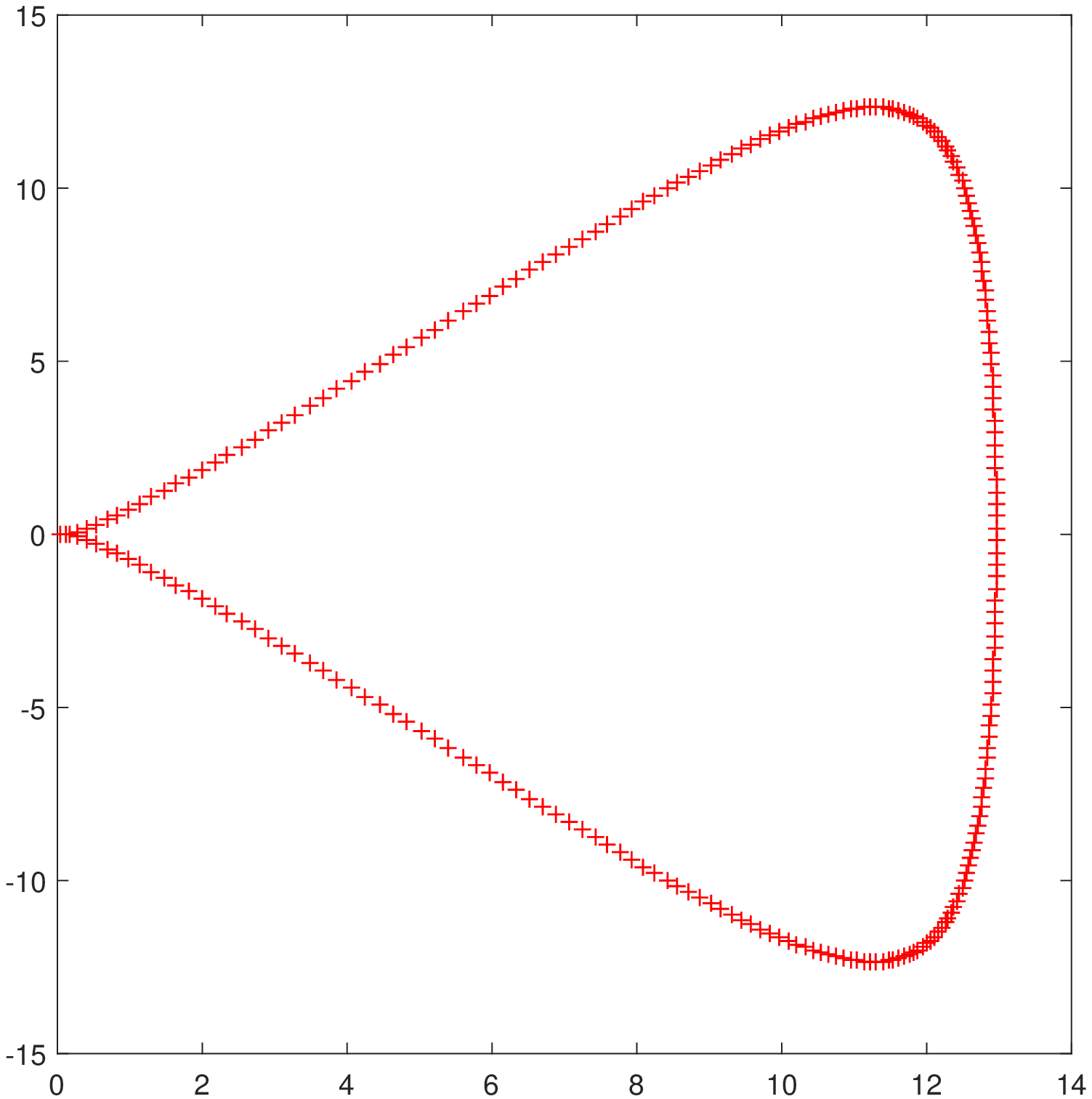}}
	\subfigure[Eigenvalues of $P_{sk}^{-1} A_0$ ({\color{blue} $*$})
	and $P_{s}^{-1} A_0$ ({\color{magenta} $\times$})]
	{\includegraphics[width=3.0in,height=2.2in]{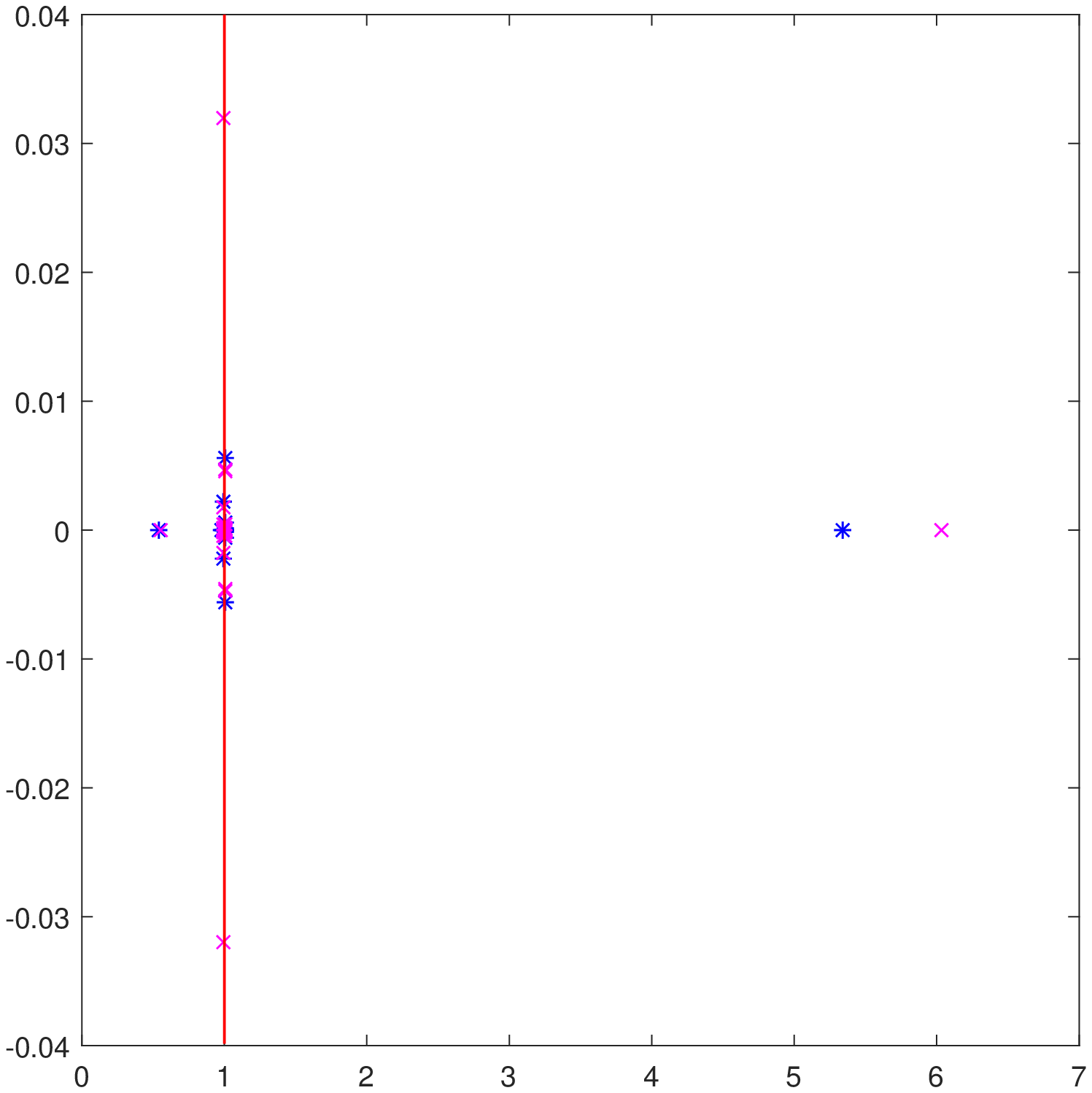}}
	\caption{Spectra of $A_0$, $P_{s}^{-1} A_0$ and $P_{s}^{-1} A_0$, when $M = N = 2^8$ in Example 1.
		Top row: $(\alpha, \beta) = (0.1, 1.1)$; Bottom row: $(\alpha, \beta) = (0.7, 1.4)$.}
	\label{fig3}
\end{figure}
\begin{table}[th]\footnotesize\tabcolsep=2.0pt
	\begin{center}
		\caption{Results of different methods when $M = N$ for Example 1.}
		\centering
		\begin{tabular}{cccccccccccc}
			\hline
			& & \rm{BS} & \rm{BFSM} & \multicolumn{2}{c}{\rm{SK2-PBiCGSTAB}}
			& \multicolumn{2}{c}{\rm{S2-PBiCGSTAB}} & \multicolumn{2}{c}{\rm{SK2-FGMRES}}
			& \multicolumn{2}{c}{\rm{S2-FGMRES}} \\
			[-2pt] \cmidrule(lr){5-6} \cmidrule(lr){7-8} \cmidrule(lr){9-10} \cmidrule(lr){11-12}\\ [-11pt]
			($\alpha$, $\beta$) & $N$ & $\mathrm{Time}$ & $\mathrm{Time}$
			&($\mathrm{Iter}$, $\mathrm{Iter3}$) & $\mathrm{Time}$
			& ($\mathrm{Iter}$, $\mathrm{Iter3}$) & $\mathrm{Time}$
			& ($\mathrm{Iter}$, $\mathrm{Iter3}$) & $\mathrm{Time}$
			& ($\mathrm{Iter}$, $\mathrm{Iter3}$) & $\mathrm{Time}$ \\
			\hline
			(0.1, 1.1) & 64 & 0.213 & 0.007 & (4+2, 5) & 0.014 & (5+2, 5) & 0.015 & (6+5, 5) & 0.020 & (6+5, 6) & 0.021 \\
			& 128 & 3.469 & 0.044 & (4+2, 5) & 0.056 & (5+2, 5) & 0.057 & (6+5, 5) & 0.077 & (6+6, 5) & 0.092 \\
			& 256 & 237.015 & 0.234 & (5+2, 5) & 0.142 & (5+2, 5) & 0.144 & (6+6, 5) & 0.234 & (6+7, 5) & 0.273 \\
			& 512 & \dag & 1.839 & (5+2, 5) & 0.995 & (5+2, 5) & 0.998 & (6+7, 5) & 1.912 & (6+8, 5) & 2.185 \\
			& 1024 & \dag & 19.839 & (5+2, 5) & 2.635 & (5+2, 6) & 2.672 & (6+9, 5) & 6.672 & (6+10, 5) & 7.480\\
			\\														
			(0.4, 1.7) & 64 & 0.185 & 0.009 & (4+2, 5) & 0.014 & (6+2, 6) & 0.015 & (6+5, 7) & 0.021 & (7+5, 6) & 0.022 \\
			& 128 & 2.993 & 0.043 & (4+2, 5) & 0.057 & (6+2, 5) & 0.058 & (6+6, 6) & 0.090 & (7+5, 6) & 0.078 \\
			& 256 & 232.214 & 0.235 & (6+2, 5) & 0.140 & (6+2, 5) & 0.141 & (6+7, 6) & 0.268 & (7+6, 5) & 0.233 \\
			& 512 & \dag & 1.840 & (6+3, 5) & 1.486 & (6+3, 5) & 1.485 & (6+7, 5) & 1.906 & (7+6, 5) & 1.664 \\
			& 1024 & \dag & 19.838 & (6+3, 5) & 3.887 & (6+3, 5) & 3.878 & (6+8, 5) & 5.983 & (7+7, 5) & 5.248 \\
			\\														
			(0.7, 1.4) & 64 & 0.183 & 0.009 & (4+3, 5) & 0.020 & (5+3, 5) & 0.020 & (6+6, 7) & 0.024 & (6+6, 8) & 0.025 \\
			& 128 & 2.969 & 0.040 & (5+3, 5) & 0.081 & (5+3, 5) & 0.083 & (6+7, 6) & 0.104 & (7+8, 6) & 0.119 \\
			& 256 & 237.030 & 0.238 & (5+4, 5) & 0.279 & (5+4, 5) & 0.279 & (6+8, 6) & 0.300 & (7+9, 6) & 0.342 \\
			& 512 & \dag & 1.842 & (5+4, 5) & 1.975 & (5+4, 5) & 1.988 & (6+10, 5) & 2.688 & (7+11, 6) & 2.971 \\
			& 1024 & \dag & 19.847 & (5+5, 5) & 6.429 & (5+5, 5) & 6.526 & (6+11, 5) & 8.174 & (7+14, 5) & 10.540 \\
			\\														
			(0.9, 1.9) & 64 & 0.176 & 0.009 & (4+2, 5) & 0.015 & (6+2, 5) & 0.016 & (5+5, 5) & 0.200 & (6+5, 5) & 0.021 \\
			& 128 & 2.950 & 0.043 & (6+3, 5) & 0.081 & (6+3, 5) & 0.082 & (6+6, 5) & 0.091 & (6+6, 5) & 0.092 \\
			& 256 & 233.143 & 0.209 & (6+3, 5) & 0.209 & (6+3, 5) & 0.214 & (6+7, 5) & 0.267 & (6+7, 5) & 0.271 \\
			& 512 & \dag & 1.837 & (6+4, 5) & 1.968 & (6+4, 5) & 1.986 & (6+8, 5) & 2.164 & (6+8, 5) & 2.182 \\
			& 1024 & \dag & 19.853 & (6+4, 5) & 5.211 & (6+4, 5) & 5.276 & (6+10, 5) & 7.505 & (6+10, 5) & 7.447 \\
			\hline
		\end{tabular}
		\label{tab1}
	\end{center}
\end{table}
\begin{table}[H]\footnotesize\tabcolsep=2.0pt
	\begin{center}
		\caption{Results of different methods when $M = 257$ for Example 1.}
		\centering
		\begin{tabular}{cccccccccccc}
			\hline
			& & \rm{BS} & \rm{BFSM} & \multicolumn{2}{c}{\rm{SK2-PBiCGSTAB}}
			& \multicolumn{2}{c}{\rm{S2-PBiCGSTAB}} & \multicolumn{2}{c}{\rm{SK2-FGMRES}}
			& \multicolumn{2}{c}{\rm{S2-FGMRES}} \\
			[-2pt] \cmidrule(lr){5-6} \cmidrule(lr){7-8} \cmidrule(lr){9-10} \cmidrule(lr){11-12}\\ [-11pt]
			($\alpha$, $\beta$) & $N$ & $\mathrm{Time}$ & $\mathrm{Time}$
			&($\mathrm{Iter}$, $\mathrm{Iter3}$) & $\mathrm{Time}$
			& ($\mathrm{Iter}$, $\mathrm{Iter3}$) & $\mathrm{Time}$
			& ($\mathrm{Iter}$, $\mathrm{Iter3}$) & $\mathrm{Time}$
			& ($\mathrm{Iter}$, $\mathrm{Iter3}$) & $\mathrm{Time}$ \\
			\hline
			(0.1, 1.1) & 65 & 3.198 & 0.077 & (4+2, 5) & 0.053 & (5+2, 5) & 0.055 & (6+5, 5) & 0.075 & (6+5, 6) & 0.072 \\
			& 129 & 13.545 & 0.115 & (4+2, 5) & 0.079 & (5+2, 5) & 0.079 & (6+5, 5) & 0.111 & (6+6, 5) & 0.138 \\
			& 257 & 277.409 & 0.209 & (5+2, 5) & 0.138 & (5+2, 5) & 0.138 & (6+6, 5) & 0.236 & (6+7, 5) & 0.269 \\
			& 513 & \dag & 0.819 & (5+2, 5) & 0.236 & (5+2, 5) & 0.237 & (6+7, 5) & 0.462 & (6+8, 5) & 0.527 \\
			& 1025 & \dag & 4.613 & (5+2, 5) & 0.405 & (5+2, 6) & 0.418 & (6+9, 5) & 1.078 & (6+10, 5) & 1.205 \\
			\\											
			(0.4, 1.7) & 65 & 3.116 & 0.072 & (4+2, 5) & 0.050 & (6+2, 6) & 0.059 & (6+5, 7) & 0.066 & (7+6, 6) & 0.077 \\
			& 129 & 13.397 & 0.119 & (4+2, 5) & 0.080 & (6+2, 5) & 0.084 & (6+6, 6) & 0.131 & (7+6, 6) & 0.137 \\
			& 257 & 263.419 & 0.210 & (6+2, 5) & 0.139 & (6+2, 5) & 0.139 & (6+7, 6) & 0.271 & (7+6, 5) & 0.238 \\
			& 513 & \dag & 0.816 & (6+2, 5) & 0.233 & (6+2, 5) & 0.232 & (6+7, 5) & 0.474 & (7+6, 5) & 0.421 \\
			& 1025 & \dag & 4.613 & (6+2, 5) & 0.410 & (6+3, 5) & 0.615 & (6+7, 5) & 0.840 & (7+6, 5) & 0.760 \\
			\\											
			(0.7, 1.4) & 65 & 3.056 & 0.073 & (4+4, 5) & 0.093 & (5+4, 5) & 0.100 & (6+8, 6) & 0.101 & (6+8, 7) & 0.109 \\
			& 129 & 13.421 & 0.115 & (5+4, 5) & 0.155 & (5+4, 5) & 0.165 & (6+8, 6) & 0.172 & (7+9, 6) & 0.199 \\
			& 257 & 251.611 & 0.214 & (5+4, 5) & 0.269 & (5+4, 5) & 0.277 & (6+8, 6) & 0.298 & (7+9, 6) & 0.334 \\
			& 513 & \dag & 0.833 & (5+4, 5) & 0.450 & (5+4, 5) & 0.457 & (6+9, 6) & 0.593 & (7+10, 6) & 0.658 \\
			& 1025 & \dag & 4.397 & (5+4, 5) & 0.792 & (5+4, 5) & 0.793 & (7+10, 5) & 1.203 & (7+11, 6) & 1.332 \\
			\\											
			(0.9, 1.9) & 65 & 3.057 & 0.070 & (4+3, 5) & 0.074 & (6+3, 5) & 0.071 & (6+7, 5) & 0.088 & (6+7, 5) & 0.102 \\
			& 129 & 13.393 & 0.118 & (4+3, 5) & 0.116 & (6+3, 5) & 0.124 & (6+7, 5) & 0.152 & (6+7, 5) & 0.160 \\
			& 257 & 257.493 & 0.211 & (6+3, 5) & 0.201 & (6+3, 5) & 0.213 & (6+7, 5) & 0.263 & (6+7, 5) & 0.274 \\
			& 513 & \dag & 0.828 & (6+3, 5) & 0.340 & (6+3, 5) & 0.365 & (6+7, 5) & 0.465 & (6+7, 5) & 0.482 \\
			& 1025 & \dag & 4.625 & (6+3, 5) & 0.587 & (6+3, 5) & 0.613 & (6+7, 5) & 0.854 & (6+7, 5) & 0.860 \\
			\hline
		\end{tabular}
		\label{tab2}
	\end{center}
\end{table}
\begin{table}[t]\footnotesize\tabcolsep=2.0pt
	\begin{center}
		\caption{Results of different methods when $N = 257$ for Example 1.}
		\centering
		\begin{tabular}{cccccccccccc}
			\hline
			& & \rm{BS} & \rm{BFSM} & \multicolumn{2}{c}{\rm{SK2-PBiCGSTAB}}
			& \multicolumn{2}{c}{\rm{S2-PBiCGSTAB}} & \multicolumn{2}{c}{\rm{SK2-FGMRES}}
			& \multicolumn{2}{c}{\rm{S2-FGMRES}} \\
			[-2pt] \cmidrule(lr){5-6} \cmidrule(lr){7-8} \cmidrule(lr){9-10} \cmidrule(lr){11-12}\\ [-11pt]
			($\alpha$, $\beta$) & $M$ & $\mathrm{Time}$ & $\mathrm{Time}$
			&($\mathrm{Iter}$, $\mathrm{Iter3}$) & $\mathrm{Time}$
			& ($\mathrm{Iter}$, $\mathrm{Iter3}$) & $\mathrm{Time}$
			& ($\mathrm{Iter}$, $\mathrm{Iter3}$) & $\mathrm{Time}$
			& ($\mathrm{Iter}$, $\mathrm{Iter3}$) & $\mathrm{Time}$ \\
			\hline
			(0.1, 1.1) & 65 & 4.065 & 0.048 & (5+2, 5) & 0.035 & (5+2, 5) & 0.037 & (6+6, 5) & 0.059 & (6+7, 5) & 0.067 \\
			& 129 & 23.278 & 0.091 & (5+2, 5) & 0.063 & (5+2, 5) & 0.070 & (6+6, 5) & 0.113 & (6+7, 5) & 0.135 \\
			& 257 & 277.409 & 0.209 & (5+2, 5) & 0.138 & (5+2, 5) & 0.138 & (6+6, 5) & 0.236 & (6+7, 5) & 0.269 \\
			& 513 & \dag & 0.569 & (5+2, 5) & 0.257 & (5+2, 5) & 0.266 & (6+6, 5) & 0.463 & (6+7, 5) & 0.550 \\
			& 1025 & \dag & 1.702 & (5+2, 5) & 0.530 & (5+2, 5) & 0.542 & (6+6, 5) & 0.958 & (6+7, 5) & 1.095 \\
			\\											
			(0.4, 1.7) & 65 & 4.057 & 0.047 & (6+2, 5) & 0.036 & (6+2, 5) & 0.039 & (6+7, 6) & 0.069 & (7+6, 5) & 0.067 \\
			& 129 & 22.929 & 0.088 & (5+2, 5) & 0.068 & (6+2, 5) & 0.071 & (6+7, 6) & 0.125 & (7+6, 5) & 0.118 \\
			& 257 & 263.419 & 0.210 & (6+2, 5) & 0.139 & (6+2, 5) & 0.139 & (6+7, 6) & 0.271 & (7+6, 5) & 0.238 \\
			& 513 & \dag & 0.585 & (6+3, 5) & 0.390 & (6+3, 5) & 0.396 & (6+7, 6) & 0.534 & (7+6, 5) & 0.476 \\
			& 1025 & \dag & 1.854 & (6+3, 5) & 0.782 & (6+3, 5) & 0.810 & (6+7, 6) & 1.101 & (7+6, 5) & 0.968 \\
			\\											
			(0.7, 1.4) & 65 & 4.069 & 0.048 & (5+3, 5) & 0.051 & (5+3, 5) & 0.550 & (6+9, 6) & 0.088 & (7+8, 6) & 0.079 \\
			& 129 & 23.217 & 0.090 & (5+3, 5) & 0.099 & (5+3, 5) & 0.106 & (6+8, 6) & 0.146 & (7+8, 6) & 0.151 \\
			& 257 & 251.611 & 0.214 & (5+4, 5) & 0.269 & (5+4, 5) & 0.277 & (6+8, 6) & 0.298 & (7+9, 6) & 0.334 \\
			& 513 & \dag & 0.578 & (4+4, 5) & 0.521 & (5+4, 5) & 0.542 & (6+9, 6) & 0.680 & (7+10, 6) & 0.776 \\
			& 1025 & \dag & 1.763 & (4+5, 5) & 1.301 & (5+5, 5) & 1.318 & (6+11, 5) & 1.680 & (7+12, 6) & 1.829 \\
			\\											
			(0.9, 1.9) & 65 & 4.081 & 0.043 & (6+2, 5) & 0.036 & (6+2, 5) & 0.039 & (6+5, 5) & 0.049 & (6+5, 5) & 0.053 \\
			& 129 & 23.093 & 0.095 & (6+3, 5) & 0.098 & (6+3, 5) & 0.106 & (6+6, 5) & 0.112 & (6+6, 5) & 0.119 \\
			& 257 & 257.493 & 0.211 & (6+3, 5) & 0.201 & (6+3, 5) & 0.213 & (6+7, 5) & 0.263 & (6+7, 5) & 0.274 \\
			& 513 & \dag & 0.552 & (6+4, 5) & 0.513 & (6+4, 5) & 0.539 & (6+8, 5) & 0.612 & (6+8, 5) & 0.642 \\
			& 1025 & \dag & 1.725 & (4+4, 5) & 1.033 & (6+4, 5) & 1.082 & (6+10, 5) & 1.528 & (6+10, 5) & 1.567 \\
			\hline
		\end{tabular}
		\label{tab3}
	\end{center}
\end{table}
\begin{table}[H]\footnotesize\tabcolsep=3.0pt
	\begin{center}
		\caption{Comparison results of SK2-PBiCGSTAB method and Huang-Lei's method for Example 1, where $M = 257$.}
		\centering
		\begin{tabular}{cccccccc}
			\hline
			& & \multicolumn{3}{c}{\rm{Huang-Lei's method}}
			& \multicolumn{3}{c}{\rm{SK2-PBiCGSTAB}} \\
			[-2pt] \cmidrule(lr){3-5} \cmidrule(lr){7-8} \cmidrule(lr){6-8} \\ [-11pt]
			($\alpha$, $\beta$) & $N$ & $\mathrm{Time}$ & $\mathrm{Error1}$ & $\mathrm{Error2}$
			& $\mathrm{Time}$ & $\mathrm{Error1}$ & $\mathrm{Error2}$\\
			\hline
			(0.1, 1.1) & 65 & 0.058 & 8.3526E-04 & 5.9916E-04 & 0.053 & 8.3526E-04 & 5.9916E-04 \\
			& 129 & 0.065 & 2.1165E-04 & 1.5173E-04 & 0.079 & 2.1165E-04 & 1.5173E-04 \\
			& 257 & 0.078 & 5.2851E-05 & 3.7902E-05 & 0.138 & 5.2852E-05 & 3.7903E-05 \\
			& 513 & 0.116 & 1.2783E-05 & 9.2066E-06 & 0.236 & 1.2778E-05 & 9.2035E-06 \\
			& 1025 & 0.233 & 2.7253E-06 & 2.0070E-06 & 0.405 & 2.7131E-06 & 1.9997E-06 \\
			\\							
			(0.4, 1.7) & 65 & 0.055 & 5.4781E-04 & 3.8003E-04 & 0.050 & 5.4781E-04 & 3.8003E-04 \\
			& 129 & 0.063 & 1.3690E-04 & 9.5128E-05 & 0.080 & 1.3689E-04 & 9.5126E-05 \\
			& 257 & 0.076 & 3.2744E-05 & 2.2885E-05 & 0.139 & 3.2743E-05 & 2.2884E-05 \\
			& 513 & 0.115 & 6.6208E-06 & 4.7452E-06 & 0.233 & 6.6207E-06 & 4.7454E-06 \\
			& 1025 & 0.224 & 1.5886E-06 & 4.9796E-07 & 0.410 & 1.5886E-06 & 4.9764E-07 \\
			\\							
			(0.7, 1.4) & 65 & 0.053 & 7.0888E-04 & 4.9767E-04 & 0.093 & 7.0888E-04 & 4.9767E-04 \\
			& 129 & 0.064 & 1.7789E-04 & 1.2502E-04 & 0.155 & 1.7790E-04 & 1.2502E-04 \\
			& 257 & 0.078 & 4.3826E-05 & 3.0074E-05 & 0.269 & 4.3825E-05 & 3.0076E-05 \\
			& 513 & 0.120 & 1.1377E-05 & 6.1321E-06 & 0.450 & 1.1376E-05 & 6.1350E-06 \\
			& 1025 & 0.220 & 2.9060E-06 & 5.7145E-07 & 0.792 & 2.9113E-06 & 5.4756E-07 \\
			\\							
			(0.9, 1.9) & 65 & 0.053 & 4.4937E-04 & 3.1623E-04 & 0.074 & 4.4937E-04 & 3.1623E-04 \\
			& 129 & 0.064 & 1.1041E-04 & 7.7685E-05 & 0.116 & 1.1043E-04 & 7.7700E-05 \\
			& 257 & 0.092 & 2.5058E-05 & 1.7763E-05 & 0.201 & 2.5028E-05 & 1.7741E-05 \\
			& 513 & 0.126 & 3.8914E-06 & 2.8666E-06 & 0.340 & 3.8553E-06 & 2.8317E-06 \\
			& 1025 & 0.220 & 1.7111E-06 & 1.0294E-06 & 0.587 & 1.7104E-06 & 1.0289E-06 \\
			\hline
		\end{tabular}
		\label{tab4}
	\end{center}
\end{table}
\begin{table}[t]\footnotesize\tabcolsep=2.0pt
	\begin{center}
		\caption{The condition numbers of $W$, $P_{W}^{-1} W$, $A_0$, $P_{s}^{-1} A_0$
			and $P_{sk}^{-1} A_0$ in Example 1.}
		\centering
		\begin{tabular}{ccccccc}
			\hline
			($\alpha$, $\beta$) & ($N$, $M$) & $W$ & $P_{W}^{-1} W$
			& $A_0$ & $P_{s}^{-1} A_0$ & $P_{sk}^{-1} A_0$ \\
			\hline
			(0.1, 1.1) & (32, 32) & 27.98 & 1.01 & 25.28 & 99.15 & 14.16 \\
			& (64, 32) & 57.43 & 1.01 & 51.90 & 212.95 & 27.82 \\
			& (128, 32) & 120.74 & 1.01 & 109.09 & 457.09 & 57.03 \\
			\\						
			(0.4, 1.7) & (32, 32) & 214.57 & 1.02 & 132.85 & 223.71 & 49.84 \\
			& (64, 32) & 696.64 & 1.02 & 431.24 & 725.02 & 152.98 \\
			& (128, 32) & 2262.94 & 1.02 & 1400.75 & 2348.38 & 484.23 \\
			\\						
			(0.7, 1.4) & (32, 32) & 89.65 & 1.05 & 39.59 & 40.06 & 18.52 \\
			& (64, 32) & 236.56 & 1.05 & 104.15 & 102.99 & 45.01 \\
			& (128, 32) & 624.16 & 1.05 & 274.49 & 268.20 & 114.37 \\
			\\						
			(0.9, 1.9) & (32, 32) & 51.45 & 1.15 & 233.76 & 211.90 & 74.67 \\
			& (64, 32) & 3063.80 & 1.02 & 872.64 & 774.19 & 259.89 \\
			& (128, 32) & 11438.08 & 1.02 & 3256.96 & 2854.62 & 932.00 \\
			\hline
		\end{tabular}
		\label{tab5}
	\end{center}
\end{table}

In Tables \ref{tab1}-\ref{tab3}, compared with BS method, the four preconditioned iterative methods
(i.e., SK2-PBiCGSTAB, SK2-FGMRES, S2-PBiCGSTAB and S2-FGMRES)
greatly reduce the computational cost in aspects of CPU time and memory requirement.
When $M = N = 2^6, 2^7$ and $2^8$ in Table \ref{tab1},
although the four preconditioned iterative methods are slower than BFSM method,
they do not need to deal with $M$ systems.
After further investigating Tables \ref{tab1}-\ref{tab3}, we have found that
there is little difference in the CPU time and number of iterations
between SK2-PBiCGSTAB and S2-PBiCGSTAB (or between SK2-FGMRES and S2-FGMRES).
However, $\textrm{Time}$ and number of iterations needed by SK2-FGMRES (or S2-FGMRES) are
slightly larger than SK2-PBiCGSTAB (or S2-PBiCGSTAB).
In Table \ref{tab4}, the SK2-PBiCGSTAB method is compared with
	the method proposed in \cite{huang2017fast} (referred to as Huang-Lei's method)
	in terms of CPU cost and accuracy of solutions.
	Here and hereafter, $\textrm{Error1} = \max\limits_{1 \leq j \leq M} \| \bm{\zeta}^j \|_\infty$ and
	$\textrm{Error2} = \max\limits_{1 \leq j \leq M} \| \bm{\zeta}^j \|$,
	where $\| \cdot \|$ is the $L_2$-norm, and $\bm{\zeta}^j$ is a vector representing the absolute error
	between the exact solution and numerical solution at $t = t_j$.
	As seen from Table \ref{tab4}, the SK2-PBiCGSTAB method needs more CPU time when solving Eq. (2.5).
	However, Error2 calculated by the SK2-PBiCGSTAB method
	is slightly smaller than the Huang-Lei's method when $N$ becomes increasingly large.
	In Table \ref{tab5}, the condition numbers of $W$, $P_{W}^{-1} W$, $A_0$, $P_{s}^{-1} A_0$
	and $P_{sk}^{-1} A_0$ are listed to further illustrate the effectiveness of $P_{W}$ and $P_{sk}$.
	It shows that both $P_{W}$ and $P_{sk}$ reduce the condition numbers greatly,
	and $P_{sk}$ performs better than $P_{s}$.
	Meanwhile, it is also interesting to notice that the condition number of $P_{s}^{-1} A_0$ is
	even larger than $A_0$ when $(\alpha, \beta) = (0.1, 1.1)$ and $(0.4, 1.7)$.
Furthermore, Fig. \ref{fig2} shows the eigenvalues of $W$ and $P_{W}^{-1} W$,
when $M = N = 2^6$ and $(\alpha, \beta) = (0.1, 1.1), (0.7, 1.4)$.
Fig. \ref{fig3} is plotted to further illustrate that $P_{sk}$
is slightly better than the Strang's preconditioner $P_{s}$.

\section{Concluding remarks}
\label{sec5}

The BLTT system (2.5) arising from TSFDE (\ref{eq1.1}) is studied.
Firstly, the $L2$-$1_\sigma$ and WSGD formulae are adopted to discrete (\ref{eq1.1}).
Secondly, for the purpose of fast solving the obtained BLTT system (2.5),
two preconditioners (i.e., $P_W$ and $P_{sk}$) are proposed and analyzed, respectively.
Finally, numerical experiments show that our proposed SK2 strategy is efficient
for fast solving the BLTT system.
Meanwhile, the numerical experiments also indicate that
the performance of our skew-circulant preconditioner $P_{sk}$ is slightly better than
	the Strang's circulant preconditioner $P_s$.
Based on this research, we give three future research directions:
(i) Notice that the preconditioner $P_W$ only compresses the temporal component.
Hence, it is valuable to develop a preconditioner which compresses both the temporal and spatial components;
(ii) $P_W$ is not suitable for parallel computing.
Thus, it is interesting to design an efficient and parallelizable preconditioner;
(iii) Some other applications of our new skew-circulant preconditioner are worth considering.

\section*{Acknowledgments}
\addcontentsline{toc}{section}{Acknowledgments}
\label{sec6}

\textit{The authors would like to thank Dr. Jiwei Zhang and Dr. Meng Li for giving some helpful discussions.
	We would like to express our sincere thanks to the referees
	and our editor Prof. Michael Ng for insightful comments and invaluable suggestions
	that greatly improved the presentation of this paper.
	We are also grateful to Dr. Siu-Long Lei for sharing us with MATLAB codes of Ref. \cite{huang2017fast}.
	This research is supported by the National Natural Science Foundation of China (Nos. 61876203, 61772003 and 11801463)
	and the Fundamental Research Funds for the Central Universities (Nos. ZYGX2016J132 and JBK1809003).
}


\end{document}